\documentclass[a4paper, 11pt]{amsart}
\usepackage{color,hyperref}
\hypersetup{
	colorlinks,
	linkcolor={red!50!black},
	citecolor={blue!50!black},
	urlcolor={blue!80!black}
}

\usepackage{graphicx}
\usepackage{xcolor}
\usepackage[margin=3cm]{geometry}
\usepackage{mathtools}
\usepackage{amsmath}
\usepackage{amssymb}
\usepackage[utf8]{inputenc}
\usepackage{graphicx}
\usepackage{stmaryrd}
\usepackage{bbold}
\usepackage[margin=3cm]{geometry}
\usepackage{float}
\usepackage[utf8]{inputenc}
\usepackage{amsthm}
\usepackage{mathtools}
\usepackage{amssymb}
\usepackage{amsmath}
\usepackage{blkarray}
\usepackage{color,hyperref}
\usepackage{enumerate}

\newcommand{\N}{\mathbb N}

\newcommand{\A}{\mathbb{A}}
\newcommand{\B}{\mathbb{B}}
\newcommand{\F}{\mathbb F}
\newcommand{\C}{\mathbb C}
\newcommand{\D}{\mathbb D}
\newcommand{\X}{\mathbb X}
\newcommand{\Y}{\mathbb Y}

\newcommand{\op}{\normalfont\text{op}}

\newcommand{\codisc}{\mathbf{indisc}}
\newcommand{\disc}{\mathbf{disc}}

\newcommand{\E}{\mathcal{E}}

\newcommand{\s}{\mathbf{Set}}

\newcommand{\GpdE}{\mathbf{Gpd(\E)}}
\newcommand{\Cat}{\mathbf{Cat}}
\newcommand{\CatE}{\mathbf{Cat(\E)}}

\newcommand{\U}{\mathcal{U}}
\newcommand{\K}{\mathcal{K}}
\newcommand{\FG}{\mathbb{F}\mathcal{G}}

\theoremstyle{plain}
\newtheorem{thm}{Theorem}[section]
\newtheorem{lem}[thm]{Lemma}
\newtheorem{cor}[thm]{Corollary}

\newtheorem{prop}[thm]{Proposition}

\theoremstyle{definition}
\newtheorem{example}[thm]{Example}
\newtheorem{define}[thm]{Definition}

\newtheorem{remark}[thm]{Remark}

\usepackage{tikz-cd}
\usepackage{tikz}
\usetikzlibrary{arrows}
\usepackage{pgfplots}
\pgfplotsset{compat=1.18}
\usetikzlibrary{positioning, calc}
\usetikzlibrary{decorations.markings,arrows.meta,bending}
\tikzset{-->--/.style={decoration={markings, % https://tex.stackexchange.com/a/39282/121799
			mark=at position #1 with {\arrow[line width=2pt]{>}}},postaction={decorate}}}
\tikzset{shorten <>/.style={shorten >=#1,shorten <=#1}}

\NewDocumentCommand\Cycle{O{} m m m O{} m}{%
	\draw[#1](#2.{#3+asin(#6/(#4*1.41))}) arc (180+#3-45:180+#3-45-270:#6/2) #5;
}

\tikzset{
	partial ellipse/.style args={#1:#2:#3}{
		insert path={+ (#1:#3) arc (#1:#2:#3)}
	}
}
\tikzset{->-/.style={decoration={
			markings,
			mark=at position .5 with {\arrow{>}}},postaction={decorate}}}

\usepackage{tkz-euclide}
\usepackage{tikz-3dplot}
\tdplotsetmaincoords{60}{60}

\usepackage{stackengine, graphicx}

\thanks{This paper was written whilst the first named author was supported by the Dame Kathleen Ollerenshaw PhD studentship and the second named author was supported by EPSRC under grant EP/V002325/2. We are both grateful for the financial support. Both authors would like to thank Nicola Gambino for his advice while preparing this paper, and to helpful comments after presenting this work at PSSL 110 in Tallinn. In particular, we would like to thank Nathanael Arkor and Peter Lefanu Lumsdaine for helpful discussions and historical insight.}

\title{Colimits of Internal Categories}
\author{Calum Hughes and Adrian Miranda}

\address{Department of Mathematics, University of Manchester, Alan Turing building, Oxford Road, Manchester M13 9PL, United Kingdom}

\email{calum.hughes77@gmail.com \\ adriantosharmiranda@hotmail.com }

\begin{document}

\begin{abstract}
  We show that for an extensive $1$-category $\mathcal{E}$ with pullback-stable coequalisers admitting free internal categories over internal graphs, the $2$-category $\mathbf{Cat}(\mathcal{E})$ of internal categories, functors and natural transformations has finite $2$-colimits. In addition, $\mathbf{Cat}(\mathcal{E})$ is extensive and codescent coequalisers are stable under pullback along discrete Conduch\'{e} fibrations. Moreover, we give converse results to this. \end{abstract}

\maketitle
\tableofcontents

\section{Introduction}

\subsection{Context and Motivation}

Whilst coproducts in $\Cat$ are easily calculated by working levelwise, the coequaliser of a parallel pair of functors $F, G: \mathcal{C} \to \mathcal{D}$ has a more complicated description involving not just equivalence classes of objects and morphisms of $\mathcal{D}$ but also equivalence classes of paths, as described in \cite{bednarczyk1999generalized}. Together, these constructions ensure that $\Cat$ has all finite colimits, which is a foundational result in the field.

The goal of this work is to provide conditions on a category $\E$ such that the $2$-category $\CatE$ of internal categories, internal functors, and internal natural transformations has finite $2$-colimits. In order to show that a $2$-category has finite $2$-colimits it suffices to show that it has coproducts, copowers by the free-living arrow in $\Cat$ (which we denote $\mathbf{2}$) and coequalisers \cite[\S 3]{kelly1989elementary}.  Extensivity of an $\E$ with pullbacks suffices for coproducts and copowers by $\mathbf{2}$ to exist in $\CatE$, as shown in Lemma 5.2 and Theorem 5.5 of \cite{hughes2024elementarytheory2categorysmall} and reviewed in Section \ref{coproducts and copowers}. We show that for a $1$-category $\E$, the property of having pullback-stable coequalisers (Definition~\ref{def: pullback stable coequalisers}) gives rise to very special coequalisers in $\CatE$, as treated in Section \ref{sec: coequalisers that agree on objects}. However, the following example illustrates that the property of having pullback-stable coequalisers in $\E$ is insufficient for $\CatE$ to have coequalisers. 

\begin{example}\label{Example coequaliser not in FinCat}
    Let $\E :=\mathbf{FinSet}$, the elementary topos of finite sets. Note that elementary toposes are locally cartesian closed, and so have pullback stable coequalisers in the sense of Definition~\ref{def: pullback stable coequalisers}. Consider the following diagram in $\CatE$, in which the two functors pick out the source and target of the free-living arrow.  
    
    \begin{equation}\label{diagram: free monoid}
        \begin{tikzcd}
        \mathbf{1} \arrow[rr, shift left = 2, "d^{1}"]\arrow[rr, shift right = 2, "d^{0}"'] && \mathbf{2}
    \end{tikzcd}
    \end{equation}

If a coequaliser to this diagram existed in $\Cat(\mathbf{FinSet})$, then it would have a single object since the $d^1$ and $d^0$ would cause the source and target of $\mathbf{2}$ to be glued together.  Given any finite monoid on $\mathbf{1}$, we can send the unique arrow of $\mathbf{2}$ to any of the generators of the monoid; this will coequalise Diagram~\ref{diagram: free monoid} and so induces a unique arrow out of the coequaliser by its universal property. Therefore the coequaliser of Diagram~\ref{diagram: free monoid}, if it exists, provides the free monoid on the terminal. Remark D5.3.4 of \cite{johnstone2002sketches} explains that in an elementary topos $\E$, the existence of such a free monoid is equivalent to $\E$ having a natural numbers object. The category $\mathbf{FinSet}$ does not have a natural numbers object; this can be deduced from the fact that the natural numbers object satisfies the axioms for Peano arithmetic \cite[Theorem 2]{lawvere2005elementary}, which cannot be satisfied by any finite set. Therefore, the coequaliser of this diagram does not exist in $\Cat(\mathbf{FinSet})$. In fact, the coequaliser of this diagram in $\Cat(\s)$ is given by the monoid of natural numbers, which is equivalently the free monoid on the singleton set.
\end{example}

The above example shows that having free monoids is important in the construction of coequalisers in $\CatE$. In the absence of cartesian closedness and a subobject classifier, having a \emph{parametrised list object} on $A \in \E$ implies the existence of the free monoid on $A$. This follows from \cite[Proposition 7.3]{maietti2010joyal} by restricting the construction of the free internal category on a free internal graph to one object categories and graphs. Since we want to work in a general setting, we will assume that free internal categories on internal graphs exist in $\E$ on top of preservation of coequalisers under pullback.

On the other hand, if we assume that $\E$ is locally finitely presentable, then the existence of $2$-colimits in $\CatE$ is relatively easy to prove.

\begin{prop}
\label{prop: CatE locally presentable}
    Let $\E$ be accessible. Then $\CatE$ is accessible as a $1$-category. Furthermore, if $\E$ also has finite colimits (so is locally finitely presentable), then $\CatE$ has $2$-colimits. 
\end{prop}

\begin{proof}
    Recall that $\CatE$ is of the form $\mathbf{Mod}(\mathcal{S}, \E),$ the category of models for a finite limit sketch $\mathcal{S}$ in $\E$. As $\E$ is accessible, we can apply \cite[Proposition 5.13]{lack2023virtual} and deduce that $\mathbf{Mod}(\mathcal{S}, \E)$ is accessible. For $\E$ locally finitely presentable, we instead apply Proposition 1.53 of \cite{rosicky1994locally}, and conclude that $\CatE_1$ is locally finitely presentable, so has finite colimits, in particular coequalisers. Therefore, $\CatE$ has finite $2$-colimits. 
\end{proof}

We restrict ourselves to the elementary setting of an extensive category $\E$ with pullbacks and pullback-stable coequalisers in which the forgetful functor $\mathcal{U}: \CatE_1 \to \mathbf{Gph}(\E)$ has a left adjoint. Our main result is Theorem \ref{thm: CatE has coequalisers}, which says that $\CatE$ has coequalisers under these assumptions. Finite $2$-colimits follow as a consequence. Moreover, the assumed properties of $\E$ imply certain properties of $\CatE$. In Section~\ref{sec: characterisation}, we investigate the necessity of the assumed properties on $\E$ by showing that certain $2$-categorical properties of $\CatE$ imply these assumptions on $\E$. This is all brought together in Theorem~\ref{thm: characterisation}, which shows an equivalence of extensive categories $\E$ with pullbacks and pullback-stable coequalisers in which the forgetful functor $\mathcal{U}: \CatE_1 \to \mathbf{Gph}(\E)$ has a left adjoint and $2$-categories satisfying a list of axioms. It should be noted that this theorem is written in purely $2$-categorical terms, without reference to the fact that $\K \simeq \CatE$, so that this gives an elementary characterisation of such $2$-categories, using the work of \cite{bourke2010codescent}.

Examples of categories $\E$ satisfying the necessary condition are given in Section~\ref{sec: prelim} and include elementary toposes with natural numbers objects. Therefore this work gives a generalisation of the work in \cite[Corollary 6.10]{johnstone1978algebraic}\footnote{We thank Nathanael Arkor for pointing this reference out to us, as well as very helpful discussions.}. In particular, we generalise the result from the setting of elementary toposes with natural numbers objects to a setting which does not need cartesian closedness, regularity, exactness or a subobject classifier. Our proof uses a different method to theirs, from which it is easier to understand the actual construction of coequalisers of internal categories. Moreover, our proof is self-contained and does not rely on a monadic functor theorem, potentially allowing us to generalise to more elementary settings. However, it is also of interest that internal categories are exactly the algebras for the free-forgetful monad on internal graphs, and so we provide a proof using the method of \cite[Corollary 6.10]{johnstone1978algebraic} as well in Appendix~\ref{appendix: alternative proof}. 

The study of $2$-categories of internal categories has been of increasing interest in recent years. \cite{bourke2014two} shows that the assignment $\E \mapsto \CatE$ is a kind of $2$-exact completion of the $1$-category $\E$. $2$-categories of internal categories are also of interest for matters relating to $2$-dimensional foundations of mathematics. In previous work \cite{hughes2024elementarytheory2categorysmall} we described the elementary theory of the $2$-category of small categories, which extends Lawvere's elementary theory of the category of sets to the higher dimensional setting. This will be extended in future work, where we will describe $2$-categories of categories, which should provide examples of elementary $2$-toposes. Although many possible definitions of elementary $2$-toposes have been given \cite{weber2007yoneda, street1980cosmoi, helfer2024internal}, it is generally agreed that $2$-toposes should have $2$-colimits. Hence, it is important to understand $2$-categories which have $2$-colimits, and our present work establishes this for $2$-categories of internal categories under appropriate assumptions on $\E$. Relatedly, our result allows for a proof that the model structure on internal categories described in \cite{everaert2005model} is cofibrantly generated and algebraic, in upcoming work \cite{hughes2024internalmodel}.

\subsection{Structure of the paper}

This work is divided into six sections. Section~\ref{coproducts and copowers} recalls the construction of coproducts and copowers by $\mathbf{2}$ in $\CatE$, and gives a more detailed outline of our strategy in constructing coequalisers. Section~\ref{sec: coequalisers that agree on objects} constructs coequalisers of parallel pairs of internal functors that agree on objects. This simple case allows us to construct coequifiers in $\CatE$. We prove that in this case coequalisers of parallel pairs of internal functors that agree on objects are stable under pullback along discrete Conduch\'{e} fibrations. Section~\ref{Coequalisers of pairs of arrows out of a discrete category} uses free internal categories on internal graphs to construct coequalisers of pairs of arrows out of a discrete category. Section~\ref{sec: coeqaulisers finish} brings together all these parts to prove that $\CatE$ has coequalisers for an arbitrary pair of parallel morphisms. Section~\ref{sec: characterisation} considers results in the converse direction. We isolate pullback stability of codescent coequalisers (Definition~\ref{def: codescent coequalisers}) along discrete Conduch\'{e} fibrations in $\CatE$ as being important as it is equivalent to pullback stability of coequalisers in $\E$ (Proposition~\ref{lem: pullback stable coequalisers that agree on objects} and Proposition~\ref{lem: pullback stable codescent coequalisers implies pullback stable coequalisers in E}). Theorem~\ref{thm: characterisation} states that the assumptions that $\E$ is an extensive category with pullbacks and pullback-stable coequalisers admitting free internal categories over internal graphs is equivalent to the assumptions that $\CatE$ is extensive, has finite $2$-colimits, pullbacks and that coequalisers that agree on objects are stable under pullback along discrete Conduch\'{e} fibrations. Theorem~\ref{thm: with Bourke} uses the work of \cite{bourke2010codescent} to state this in elementary terms for a general $2$-category $\K$ without reference to the fact that $\K\simeq \CatE$. We conclude in Section~\ref{sec: prelim} with examples of when our results can be applied. This involves reproving a theorem of Maietti about constructing a left adjoint to the forgetful functor $\mathcal{U}: \CatE \to \mathbf{Gph}(\E)$ in a more general setting, though the proof is the same. We also show, due to an argument by Peter LeFanu Lumsdaine, that list-arithmetic pretoposes have pullback-stable coequalisers.

\subsection{Notational and terminology}\label{sec: notation and terminology colimits}

We adopt the notation and terminology for internal categories that was established in \cite[\S 2]{hughes2024elementarytheory2categorysmall}. In particular, note that we sometimes distinguish between the $2$-category $\CatE$ of internal categories, internal functors and internal natural transformations, and its underlying $1$-category, which we denote $\CatE_1.$ For a parallel pair of arrows $f,g: X \to Y$ and an arrow $h:Y \to Z$, we say $h$ \emph{coequalises} $f$ and $g$ for the situation in which $hg=hf$. We say that $h$ is \emph{the coequaliser} of $f$ and $g$ in the universal such case. 

Important to our proof is that coequalisers are stable under pullback.

\begin{define}\label{def: pullback stable coequalisers}
    Let $\E$ be a category with pullbacks. We say that $\E$ has \emph{pullback-stable coequalisers} if it has coequalisers and for any morphism $f: X \to Y$ in $\E$ the pullback functor $f^{*}: \E/Y \to \E/X$ preserves coequalisers.
\end{define}

More explicitly, this means that if we have $f: X \to Y$ and a coequaliser diagram:

\begin{equation}\label{equation: coequaliser}
    \begin{tikzcd}
        A \arrow[r, shift left =2, "g"] \arrow[r, shift right = 2, "h"'] & B \arrow[r, "q"] & C
    \end{tikzcd}
\end{equation}

 in $\E/Y$, then the diagram :

\begin{equation}\label{diagram: pulled-back coequaliser}
    \begin{tikzcd}
        f^*(A) \arrow[r, shift left =2, " f^*(g)"] \arrow[r, shift right = 2, " f^*(h)"'] &  f^*(B) \arrow[r, " f^*(q)"] &  f^*(C)
    \end{tikzcd}
\end{equation}

is a coequaliser diagram in $\E/X$ (and hence also in $\E$).

Recall that a \emph{regular epimorphism} is a morphism that is the coequaliser of some diagram, and a \emph{regular category} is a category with finite limits, coequalisers of kernel pairs and the property that the pullback of a regular epimorphism is a regular epimorphism. 

If $\E$ has pullback stable coequalisers, then regular epimorphisms are stable under pullback. Therefore, a category with pullback stable coequalisers and equalisers is a regular category, since we are assuming that we have pullbacks and coequalisers.

On the other hand, if $\E$ is a regular category with coequalisers, this does not guarantee that we have pullback-stable coequalisers. Given a coequaliser diagram as in Equation~\ref{equation: coequaliser}, $q: B \to C$ is a regular epimorphism, so by regularity of $\E$ it follows that $f^*(q): f^*(B) \to f^*(C)$ is a regular epimorphism, and so it is the coequaliser of some diagram:

\begin{equation*}
    \begin{tikzcd}
        E \arrow[r, shift left =2, "k"] \arrow[r, shift right = 2, " l"'] &  f^*(B) \arrow[r, " f^*(q)"] &  f^*(C).
    \end{tikzcd}
\end{equation*}

However, there is nothing to guarantee that $E = f^*(A)$ and $k= f^*(g), l=f^*(h)$, and so this does not imply pullback stability of coequalisers in the sense of Definition~\ref{def: pullback stable coequalisers}. In our explicit description of coequalisers, it is important that we have this extra control over the pullback of coequalisers.

\section{Constructing finite $2$-colimits of internal categories via simpler colimits}
\label{coproducts and copowers}

Recall (for example from \cite[\S 3]{kelly1989elementary}) that finite $2$-colimits can be constructed using finite coproducts, coequalisers of parallel pairs, and copowers by $\mathbf{2}$. We briefly review the construction of finite coproducts and copowers by $\mathbf{2}$ in the $2$-category $\CatE$ under the assumption that $\E$ is extensive and has pullbacks. We then outline the construction of coequalisers of parallel pairs in $\CatE$, which we will develop over the subsequent sections.

First, we describe an internal free-living arrow in $\CatE$, which we denote $\mathbf{2}_{\E}$. For any object $\A \in \CatE$, the cartesian product $\mathbf{2}_{\E} \times \A$ will have the universal property of the copower of $\A$ by $\mathbf{2}.$ The internal category $\mathbf{2}_{\E}$ can be concretely described as a truncated simplicial object, with $n$-simplices given by the $(n+2)$-fold coproduct of the terminal object $\mathbf{1} \in \E$; see Example 2.3.2 of \cite{miranda2022internal} for further details. Abstractly, it is the image of $\mathbf{2}$ under $\Cat(F): \Cat(\mathbf{FinSet}) \to \CatE$, where $F: \mathbf{FinSet} \to \E$ is the unique coproduct and terminal object preserving functor, which is described in Definition 5.4 of \cite{hughes2024elementarytheory2categorysmall}. We note that with the additional assumption of cartesian closedness, Proposition~\ref{proposition coproducts and copwers by 2 in Cat E} (2) is Theorem 5.5 (2) of \cite{hughes2024elementarytheory2categorysmall}, but this proof is more general as we only assume extensivity.

\begin{prop}\label{proposition coproducts and copwers by 2 in Cat E}
    Let $\E$ be an extensive category with pullbacks. Then $\CatE$ has \begin{enumerate}
        \item extensive coproducts which are created by $N: \CatE_{1} \to [\Delta_{\leq 3}^\text{op}, \E]$.
        \item copowers by $\mathbf{2}$, which for an internal category $\A$ are given by $\mathbf{2}_{\E}\times \A$.
    \end{enumerate}
\end{prop}

\begin{proof}

For part (1), the coproduct of a pair of internal categories $\A$ and $\B$ is given levelwise by $n\mapsto \A_{n}+ \B_{n}.$ We refer the reader to \cite[Lemma 5.2]{hughes2024elementarytheory2categorysmall} for a full proof and details. For part (2), the internal functor $\mathbf{2}_{\E} \times \A \to \B$ corresponding to an internal natural transformation $\overline{\alpha}: f \Rightarrow g: \A \to \B$ is given via the description of $\mathbf{2}_{\E}$ by two morphisms $(f_{0}, g_{0}): A_{0}+A_0 \to B_0$ and $(f_{1}, m.\alpha, g_{1}): A_1 +A_1 + A_1 \to B_1$ in $\E$. Further details can be found in \cite{miranda2022internal}.
    
\end{proof}

In light of Proposition \ref{proposition coproducts and copwers by 2 in Cat E}, to show that $\CatE$ has finite $2$-colimits it suffices to show that the $2$-category $\CatE$ has coequalisers of parallel pairs. Moreover, since $\CatE$ has powers by $\mathbf{2}$, it suffices to show that the underlying category $\CatE_{1}$ has coequalisers of parallel pairs. 

A naive attempt at constructing a coequaliser of a pair of internal functors would be to do this levelwise. We have already seen in Example \ref{Example coequaliser not in FinCat} that this does not work even internal to $\s$ since pairs of morphisms may become newly composable once a coequaliser is also taken at the level of objects. In Example \ref{Example coequaliser not in FinCat}, the single non-identity morphism of the free living arrow becomes composable with itself after gluing together its source and target; this new composite is not created by coequalising on morphisms, and so one must take the free category on the graph obtained by coequalising on objects and then morphisms.

Our construction of coequalisers of arbitrary parallel pairs of internal functors $F,~G:~\A~\to~\B$ decomposes into the following two steps.

\begin{enumerate}
    \item First restrict $F$ and $G$ along $\varepsilon_{\A}:\disc(A_0) \to \A$ and form the coequaliser $K: \B \to \D$ of the parallel pair $F\cdot\varepsilon_{\A}$ and $G\cdot\varepsilon_{\A}$. 
    
       \begin{equation*}
        \begin{tikzcd}
            \disc(A_0) \arrow[r, bend left, "F \cdot \varepsilon_{\A}", pos=0.4] \arrow[r, bend right, "G \cdot \varepsilon_{\A}"', pos=0.4] & \B \arrow[r, "K"] & \D.
        \end{tikzcd}
    \end{equation*}
    
    In Proposition \ref{coequaliser out of discrete category} we show that if $\E$ is a category with pullbacks and pullback-stable coequalisers in which the forgetful functor $\mathcal{U}: \CatE \to \mathbf{Gph}(\E)$ has a left adjoint, then coequalisers of parallel pairs of internal functors out of discrete categories exist in $\CatE$.

    \item Next, form the coequaliser $P: \D \to \C$ of the parallel pair of internal functors $KF$ and $KG$. 
    
    \begin{equation*}
        \begin{tikzcd}
            \A \arrow[r, bend left, "K\cdot F"] \arrow[r, bend right, "K\cdot G"'] & \D\arrow[r, "P"] & \C.
        \end{tikzcd}
    \end{equation*}
    
    Note that since $K$ coequalises $F.\varepsilon_\A$ and $G.\varepsilon_\A$, the functors $KF$ and $KG$ agree on objects. In Proposition \ref{prop coequaliser when agree on objects} we show that if $\E$ has pullback-stable coequalisers then $\CatE$ has coequalisers of parallel pairs of internal functors that agree on objects.
\end{enumerate}

Finally, in Section \ref{sec: coeqaulisers finish} we show that for abstract reasons these steps combine in such a way that $Q:= PK: \B \to \C$ is the coequaliser of the original parallel pair $F, G: \A \to \B$. We prove Proposition \ref{coequaliser out of discrete category}, as required for step (1) above, using the following two auxiliary constructions.

\begin{enumerate}[i]
    \item The construction of free categories on graphs. We use their universal property.
    \item The construction of coequifiers of parallel pairs of internal natural transformations. We show in Corollary \ref{Corollary coequifiers} that when $\E$ has pullback-stable coequalisers then $\CatE$ has coequifiers of arbitrary pairs of internal natural transformations.
\end{enumerate}

In step (1) above, we first forget about any morphisms in $\A$ and instead generate the coequaliser on objects and  consider the graph $\mathcal{G}$ which has equivalence classes of objects in $\B$ as objects and morphisms in $\B$ as edges. The free category on this graph gives us a category whose morphisms are strings of morphisms in $\B$ that become composable once we coequalise on objects. We require an internal functor $\B \to \mathbf{F}(\mathcal{G})$, but the construction so far only guarantees us a morphism of their underlying graphs. The final two coequifiers extend this to a morphism of graphs which respects identities and composition. 

Step (2) then considers the morphisms of $\A$, and takes the coequaliser just on morphisms. This requires only exactness properties in $\E$.

\begin{remark}
    It is interesting to compare this construction with the method used in \S 4 of \cite{bednarczyk1999generalized} in the context of $\Cat$. Let $F,G : \mathcal{A} \to \mathcal{B}.$ The construction of a coequaliser in \cite{bednarczyk1999generalized} first constructs a relation ${_F{=}_G}$ on $\mathcal{B}$ generated by $F$ and $G$ defined on objects by $a\mathbin{_F{=}_G}a  \in \mathcal{A}_0$ iff $F(a)=G(a)$ and on morphisms by $f\mathbin{_F{=}_G} f$ iff $F(f)=G(f)$. It then constructs the \emph{generalised congruence} ${_F{\simeq}_G}$ generated by this relation, which closes this relation on morphisms under some axioms. It then quotients $\mathcal{B}$ by this generalised congruence, and the result is the coequaliser. In contrast, Step (1) of our construction constructs a category in which the \emph{generalised} congruence on $\mathcal{B}$ is simply an ordinary congruence (in the standard sense of \cite{mac2013categories}, for example) on this new category. In other words, the category constructed by Step (1) is the setting in which the generalised congruence is defined. In internal category theory, one must be very careful to state precisely where things are defined. Step (2) takes the usual quotient of a category by a congruence.  

    We do not, however, attempt to define the notion of a generalised congruence on an internal category. 
\end{remark}

\section{Coequalisers of arrows that agree on objects}
\label{sec: coequalisers that agree on objects}

Throughout this section, $\E$ will be assumed to be a category with pullbacks and pullback-stable coequalisers. The goal of this section is to show that under these assumptions, the $2$-category $\CatE$ has coequalisers of pairs of internal functors $F, G: \A \to \B$ which agree on objects in the sense that the morphisms $F_{0}, G_{0}: A_{0} \to B_{0}$ are equal in $\E$. As a corollary, we find that $\CatE$ also has coequifiers under these assumptions. Finally, we show that these coequalisers are stable under pullback along discrete Conduch\'{e} fibrations.

\begin{prop}\label{prop coequaliser when agree on objects}
     Let $\E$ be a category with pullbacks and pullback-stable coequalisers. Any pair $F, G: \A \to \B$ of internal functors that agree on objects has a coequaliser in $\CatE$.
\end{prop}

\begin{proof}

Consider the limit $L$ in $\E$ of the following diagram:

\begin{equation*}
    \begin{tikzcd}
    & & L \arrow[drr, "\pi_2", dashed] \arrow[d, "\pi_1", dashed] \arrow[dll, "\pi_0"', dashed] & &  \\
        B_1 \arrow[dr, "d_0"] & & A_1 \arrow[dl, "F_0d_1"] \arrow[dr, "F_0d_0"] & & B_1 \arrow[dl, "d_1"] \\
        & B_0 & & B_0
    \end{tikzcd}
\end{equation*}

and define $\widetilde{F} : = (\pi_0{,} F_1\pi_1 {,} \pi_2), \widetilde{G}:= (\pi_0{,} G_1\pi_1 {,} \pi_2) : L \to B_3$. We define $Q_1: B_1 \to C_1$ as the coequaliser of the following parallel pair of arrows in $\E$:

\begin{equation*}
    \begin{tikzcd}[column sep = huge]
        L \arrow[r, shift left = 2, "\widetilde{F}"]\arrow[r, shift right = 2, "\widetilde{G}"'] & B_3 \arrow[r, "m^2"] & B_1.
    \end{tikzcd}
\end{equation*}

We show that $(B_0, C_1)$ turns out to be the objects of objects and morphisms for an internal category which has the universal property of the desired coequaliser. We define source and target $d_0, d_1: C_1 \to B_0$ using the universal property of the coequaliser given the commutativity of the following diagram for $(i,j)\in\{(0,2),(1,0)\}$.

\begin{equation*}
    \begin{tikzcd}[column sep = huge]
       L \arrow[r, "\widetilde{F}"] \arrow[d, "\widetilde{G}"'] \arrow[dr, "\pi_j"] & B_3 \arrow[r, "m^2"] \arrow[d, "\pi_j"] & B_1 \arrow[dd, "d_i"] \\
       B_3 \arrow[r, "\pi_j"] \arrow[d,"m^2"] & B_1 \arrow[dr, "d_i"] & \\
       B_1 \arrow[rr, "d_i"'] & & B_0
    \end{tikzcd}
\end{equation*}

We define $i: B_0 \to C_1$ as the composite

\begin{equation*}
    \begin{tikzcd}
         B_0 \arrow[r, "i"] & B_1 \arrow[r, "Q_1"] & C_1.
    \end{tikzcd}
\end{equation*}

Next, define $C_2$ as the pullback of $d_0, d_1: C_1 \to B_0$ and define $Q_{2}: B_{2} \to C_{2}$ to be induced by the universal property of the pullback, given the morphisms $Q_{1}\cdot \pi_{0}$ and $Q_{1}\cdot \pi_{1}$. Since coequalisers are stable under pullback in $\E$, it follows that the following is a coequaliser diagram.

\begin{equation*}
    \begin{tikzcd}[column sep = huge]
        L \times_{B_0} L \arrow[r, shift left = 2, "m^2 \cdot \widetilde{F} \times_{B_0} m^2 \cdot \widetilde{F}"]\arrow[r, shift right = 2, "m^2 \cdot \widetilde{G}\times m^2 \cdot \widetilde{G}"'] & B_2 \arrow[r, "Q_2"] & C_2.
    \end{tikzcd}
\end{equation*}

By the same reasoning, the following is also a coequaliser diagram. 

\begin{equation*}
    \begin{tikzcd}[column sep = huge]
        B_3 \times_{B_0} L \arrow[r, shift left = 2, "B_3 \times_{B_0} m^2 \cdot \widetilde{F}"]\arrow[r, shift right = 2, "B_3 \times m^2 \cdot \widetilde{G}"'] & B_2 \arrow[r, "B_3 \times_{B_0} Q_1"] & B_3 \times_{B_0} C_1.
    \end{tikzcd}
\end{equation*}

 In the following, the morphism $Q_{1}.m^2: B_{3}\times_{B_{0}}B_{1} \to C_{1}$ coequalises the pair $B_3 \times_{B_0} m^2\cdot \widetilde{F}, B_3 \times_{B_0} m^2 \cdot \widetilde{G}$ by definition of the maps involved and associativity of composition in $\B$. This induces an arrow $u: B_3 \times_{B_0}C_1 \to C_1$ by the universal property of the coequaliser.

\begin{equation*}
    \begin{tikzcd}[column sep = huge]
        B_3 \times_{B_0} L \arrow[r, shift left = 2, "B_3 \times_{B_0}\widetilde{F}"]\arrow[r, shift right = 2, "B_3 \times_{B_0}\widetilde{G}"'] \arrow[d, "m^3\times_{B_0} A_1 \times_{B_0} B_1"'] & B_3 \times_{B_0} B_3 \arrow[r, "B_3 \times_{B_0} m^2"] \arrow[d, "m^3"] & B_3 \times_{B_0} B_1\arrow[d, "m^2"] & \\
        L \arrow[r, shift left = 2, "\widetilde{F}"]\arrow[r, shift right = 2, "\widetilde{G}"'] & B_3 \arrow[r, "m^2"] &  B_1 \arrow[r, "Q_1"] & C_1
    \end{tikzcd}
\end{equation*}

 By definition of the coequalisers involved and associativity, the following diagram commutes.

\begin{equation*}
    \begin{tikzcd}
        L \times_{B_0} L \arrow[r, "\widetilde{F} \times_{B_0} L"] \arrow[d, "\widetilde{G}\times_{B_0} L"'] & B_3 \times_{B_0} L \arrow[d, "B_3 \times_{B_0} \widetilde{G}"] \arrow[rr, "B_3 \times_{B_0} \widetilde{F}"] & & B_3 \times_{B_0} B_3 \arrow[d, "B_3 \times_{B_0} m^2"'] \arrow[r, "m^4"] &    B_2 \arrow[d, "m"] \\
        B_3 \times_{B_0} L \arrow[r, "B_3 \times_{B_0} \widetilde{F}"] \arrow[ddd, "B_3 \times_{B_0} \widetilde{G}"']& B_3 \times_{B_0} B_3 \arrow[dr, "B_3 \times m^2"] &  &B_3 \times_{B_0} B_1  \arrow[dd, "B_3 \times_{B_0} Q_1"] \arrow[r, "m^3"] &   B_1 \arrow[ddd, "Q_1"] \\
        & &  B_3 \times_{B_0} B_1  \arrow[dr, "B_3 \times_{B_0} Q_1"] && \\
        & B_3 \times_{B_0} B_1 \arrow[dr, "m^3"] \arrow[rr, "B_3 \times_{B_0} Q_1"] &  & B_3 \times_{B_0} C_1 \arrow[dr, "u"] & \\
        B_3 \times_{B_0} B_3 \arrow[r, "m^2 \times_{B_0} m^2"'] \arrow[ur, "B_3 \times_{B_0} m^2"]& B_2 \arrow[r, "m"'] & B_1 \arrow[rr, "Q_1"'] & & C_1
    \end{tikzcd}
\end{equation*}

Therefore, we obtain a unique arrow $m: C_2 \to C_1$ such that $mQ_2 =Q_1 m$.

We claim that $\mathbb{C}:= (B_0, C_1, d_0, d_1, i, m)$ forms an internal category. The laws specifying the source and target of identity morphisms are satisfied as shown below:

\begin{equation*}
\begin{tikzcd}
    B_0 \arrow[r, "i"] \arrow[drr, "{B_0}"'] & B_1 \arrow[r, "Q_1"]\arrow[dr, "d_i"] & C_1 \arrow[d, "d_i"] \\
    & & B_0
\end{tikzcd}
        \qquad i \in \{0,1\}.
\end{equation*}

To show that the laws specifying the source and target of composite morphisms are satisfied, we appeal to the universal property of $Q_{2}$ as the coequaliser of $F_{2}$ and $G_{2}$. We show that, for $i \in \{0,1\}$, the maps $Q_2d_i m, Q_2d_i\pi_i : B_2 \to B_0$ are equal in the diagram below. Both maps clearly coequalise $F_2, G_2 : A_2 \to B_2$. By uniqueness aspect of the universal property, it follows that $d_im = d_i\pi_i$.

\begin{equation*}
    \begin{tikzcd}
        B_2 \arrow[r, "Q_2"] \arrow[ddd, "Q_2"'] \arrow[dr, "m" ] \arrow[ddr, "\pi_i"'] & C_2 \arrow[r, "m"] & C_1 \arrow[ddd, "d_i"] \\
        & B_1 \arrow[ur, "Q_1"] \arrow[ddr, "d_i"] & \\
        & B_1 \arrow[dr, "d_i"'] \arrow[d, "Q_1"']& \\
        C_2 \arrow[r, "\pi_i"'] & C_1 \arrow[r, "d_i"'] & C_2 
    \end{tikzcd}
\end{equation*}

The other axioms follow similarly; for example, the left unit law follows from the fact that by the assumption that coequalisers are closed under pullbacks, the following diagram is a coequaliser diagram:

\begin{equation*}
    \begin{tikzcd}[column sep = huge]
         B_0 \times_{B_0} L \arrow[r, shift left = 2, "B_0 \times_{B_0} m^2\cdot \widetilde{F}"] \arrow[r, shift right = 2, "B_0 \times_{B_0} m^2 \cdot \widetilde{G}"']  & B_0 \times_{B_0} B_1 \arrow[r, "B_0 \times_{B_0} Q_1"] & B_0 \times_{B_0} C_1
    \end{tikzcd}
\end{equation*}

and so we can check the left unit law by showing that the maps $$ m\cdot(i \times_{
B_0} C_1)\cdot(B_0 \times_{B_0} Q_1), \pi_1\cdot(B_0 \times_{B_0} Q_1): B_0 \times_{B_0} B_1 \to C_1$$ are equal, and since both maps clearly coequalise the diagram above, by uniqueness of the universal property, it follows that $m\cdot (i\times_{B_0} C_1) = \pi_1$.

The right unit law and associativity of composition follows using the same method; the details for associativity can be found in appendix~\ref{appendix: proof}. 

This shows that $\mathbb{C}$ is an internal category.

By definition of $d_0, d_1: C_1 \to B_0$, $i: B_0 \to C_1$ and $m: C_2 \to C_1$, it also follows that $Q:= (\text{id}_{B_0}, Q_1)$ is well-defined as an internal functor $\B \to \C$. We now show that it has the universal property of the coequaliser of $F$ and $G$.

Given 
\begin{equation*}
    \begin{tikzcd}
        \A \arrow[r, shift left = 2, "F"] \arrow[r, shift right = 2, "G"'] & \B \arrow[r, "Q"] \arrow[dr, "R"'] & \mathbb{C} \\
        & & \mathbb{D}
    \end{tikzcd}
\end{equation*}

where $RF = RG$, we define $K_1: C_1 \to D_1$ by the universal property of $C_1$ as a coequaliser given the following equalities.

\begin{align*}
    R_1 \cdot m^2 \cdot (B_1 \times_{B_0} F_1 \times_{B_0} B_1) & = m^2 \cdot R_3 \cdot (B_1 \times_{B_0} F_1 \times_{B_0} B_1)\\
    & = m^2 (R \times_{D_0} F_1 \cdot R \times_{D_0} R) \\
    & = m^2 (R \times_{D_0} G_1 \cdot R \times_{D_0} R) \\ 
    & = m^2 \cdot R_3 \cdot (B_1 \times_{B_0} G_1 \times_{B_0} B_1) \\
    & = R_1 \cdot m^2 \cdot (B_1 \times_{B_0} G_1 \times_{B_0} B_1). 
 \end{align*}
 
 The following diagrams show that $K : = (R_0, K_1)$ assembles into a functor $\mathbb{C} \to \mathbb{D}$. The diagrams make use of the universal property of $Q_1$ and $Q_2$ as coequalisers. Uniqueness of $K: \C \to \D$ follows from uniqueness of $K_1$.

\begin{equation*}
    \begin{tikzcd}
        B_0 \arrow[rr, "R_0"] \arrow[dd, "i"] \arrow[dr, "i"] & & D_0 \arrow[dd, "i"'] \\ 
        & B_1 \arrow[dr, "R_1"] \arrow[dl, "Q_1"] &\\
        C_1 \arrow[rr, "K_1"'] & & D_1
    \end{tikzcd}
    \qquad
    \begin{tikzcd}
        B_1 \arrow[r, "Q_1"] \arrow[d, "Q_1"'] \arrow[rr, bend right, "R_1"] \arrow[dr, "d_i"'] & C_1 \arrow[r, "K_1"] & D_1 \arrow[d, "d_i"] \\
        C_1 \arrow[r, "d_i"'] & B_0 \arrow[r, "R_0"'] & D_0
    \end{tikzcd} \qquad 
    \begin{tikzcd}
        B_2 \arrow[r, "Q_2"] \arrow[dd, "Q_2"]\arrow[rr,bend right,  "R_2"] \arrow[dr, "m"'] & C_2 \arrow[r, "K_2"] & D_2 \arrow[dd, "m"] \\
        & B_1 \arrow[d, "Q_1"] \arrow[dr, "R_1"] & \\
        C_2 \arrow[r, "m"] & C_1 \arrow[r, "K_1"] & D_1
    \end{tikzcd}
\end{equation*}

\color{black}
\end{proof}
 
 We will use the following result to show that coequifiers exist in $\CatE$ when $\E$ satisfies the assumptions of Proposition \ref{prop coequaliser when agree on objects} and is moreover extensive. 

\begin{lem}\label{lem:equifier equaliser}
    Let $\K$ be a $2$-category with powers by $\mathbf{2}$. Then the equifier of a parallel pair of $2$-cells $\overline{\alpha},\overline{\beta}:~f\Rightarrow~g: A \to B$ exists if and only if the equaliser of the corresponding morphisms $\hat{\alpha}, \hat{\beta} : A \to B^{\mathbf{2}}$ exists. In this case, the limits agree. 
\end{lem}
\begin{proof}
We can check this representably in $\Cat.$ Recall that an equaliser of $\hat{\alpha}, \hat{\beta} : \mathcal{A} \to \mathcal{B}^{\mathbf{2}}$ in $\Cat$ is given by the full subcategory of those $a \in \mathcal{A}$ such that $\hat{\alpha}(a)=\hat{\beta}(a)$. Similarly, recall that the equifier of $\overline{\alpha}, \overline{\beta}: f \Rightarrow g$ in $\Cat$ is given by the full subcategory of $a \in \mathcal{A}$ such that $\overline{\alpha}_a = \overline{\beta}_a.$ By definition, $\hat{\alpha}(a)=\alpha_a$ and $\hat{\beta}(a) = \beta_a$, so these define the same thing.  
\end{proof}

The corollary to follow records conditions under which the $2$-category $\CatE$ has coequifiers. These will be used in the construction of coequalisers of parallel pairs of internal functors whose domains are discrete in Section \ref{Coequalisers of pairs of arrows out of a discrete category}.

\begin{cor}\label{Corollary coequifiers}
   Let $\E$ be an extensive category with pullbacks and pullback-stable coequalisers. The $2$-category $\CatE$ has coequifiers.
\end{cor}

\begin{proof}
    Consider the parallel pair of internal natural transformations displayed below left. By Lemma~\ref{lem:equifier equaliser} applied to $\K = \CatE^{\op},$ these correspond to the parallel pair of internal functors displayed below right. Observe that both functors are given on objects by the morphism $(F_{0}, G_{0}): A_0 + A_0 \to B_{0}$. Hence the result follows from Proposition~\autoref{prop coequaliser when agree on objects}.

    $$\begin{tikzcd}
        \A \arrow[rr, shift left = 2, "F"name=A, bend left]\arrow[rr, shift right = 2,"G"'name=B, bend right] && \B
        \arrow[from=A, to=B, Rightarrow, shift right = 2, "\overline{\alpha}"', shorten = 5]\arrow[from=A, to=B, Rightarrow, shift left = 2, "\overline{\beta}", shorten = 5] &{}
    \end{tikzcd}\begin{tikzcd}
        \mathbf{2}_{\E} \times \A \arrow[rr, shift left = 2, "\hat{\alpha}", bend left]\arrow[rr, shift right = 2, "\hat{\beta}"', bend right]&& \B
    \end{tikzcd}$$
\end{proof}

\begin{remark}
    We also note that under the assumptions that $\E$ is a pretopos, $\CatE$ also has cocomma objects which are constructed in a similar way. Given a span of functors \begin{tikzcd}
        \A & \B \arrow[l, "F"'] \arrow[r, "G"] & \C,
    \end{tikzcd} their cocomma has object of objects given by $B_{0}+C_{0}$ and object of morphisms constructed using limits and coequalisers in $\E$. Specifically, first construct the limit $L$ of the diagram displayed below.

   $$ \begin{tikzcd}
        B_{1}\arrow[rd, "d_{0}"'] && A_{0}\arrow[ld, "F_{0}"] \arrow[rd, "G_{0}"']&& C_{1}\arrow[ld, "d_{1}"]
        \\
        & B_{0} && C_{0}
    \end{tikzcd}$$

    When $\E = \s$ this limit consists of a morphism $f$ in $\B$, a morphism $g$ in $\C$ and a `heteromorphism' from the target $Z$ of $f$ to the source $Y$ of $g$ whenever there is an object $X$ in $\A$ satisfying $FX = Z$ and $GX = Y$. This heteromorphism will correspond to the component on $X$ of the natural transformation forming part of the cocomma cocone. To ensure that these heteromorphisms form a natural transformation, we next form the coequaliser of a parallel pair of maps from $b, c: A_{1} \to L$. These maps are induced by the universal property of $L$, given the data displayed below left for $b$ and below right for $c$.

$$ \begin{tikzcd}[column sep = 15]
    &&A_{1}\arrow[dd, "d_{0}"]\arrow[lldd, "F_{1}"']\arrow[r, "d_{0}"]& A_{0}\arrow[dr, "G_{0}"]
    \\
    &&&&C_{0}\arrow[d, "i"]
    \\
        B_{1}\arrow[rd, "d_{0}"'] && A_{0}\arrow[ld, "F_{0}"] \arrow[rd, "G_{0}"']&& C_{1}\arrow[ld, "d_{1}"]&{}
        \\
        & B_{0} && C_{0}
    \end{tikzcd}\begin{tikzcd}[column sep = 15]
    &A_{0}\arrow[dl, "F_{0}"']
    &A_{1}\arrow[dd, "d_{0}"]\arrow[rrdd, "G_{1}"]\arrow[l, "d_{0}"']
    \\
    C_{0}\arrow[d, "i"']
    \\
        B_{1}\arrow[rd, "d_{0}"'] && A_{0}\arrow[ld, "F_{0}"] \arrow[rd, "G_{0}"']&& C_{1}\arrow[ld, "d_{1}"]
        \\
        & B_{0} && C_{0}
    \end{tikzcd}$$

    We leave details of the proof that this gives a well-defined internal category which has the universal property of a cocomma to the interested reader. Cocommas in $\CatE$ will not be needed in this paper.
\end{remark}

It is not true that coequalisers in $\E$ being stable under pullbacks implies that all coequalisers in $\CatE$ are stable under pullback; a  counterexample is given in the case that $\E =\s$ by Shulman \cite{ncatlab:exponentials}. However, from our construction it is not too hard to see that given pullback stability of coequalisers in $\E$, a certain class of coequalisers is stable under pullback along discrete Conduch\'{e} fibrations. Discrete Conduch\'{e} fibrations can be defined representably in any $2$-category, but we recall an equivalent description for internal categories below.

\begin{define}
    Let $\E$ be a category with pullbacks. An internal functor $F: \X \to \Y$ is called a  \emph{discrete Conduch\'{e} fibration} if the square displayed below is a pullback.

    \begin{equation*}
        \begin{tikzcd}
            X_2 \arrow[d, "m"'] \arrow[r, "F_2"] & Y_2 \arrow[d, "m"] \\
            X_1 \arrow[r, "F_1"'] & Y_1
        \end{tikzcd}
    \end{equation*}
\end{define}

Discrete Conduch\'{e} fibrations are precisely those internal functors in which we can lift composites in $\Y$ that are in the image of $F$ uniquely to composites in $\X$. The following result will be useful.

\begin{lem}\label{lem: conduche}
    Let $F: \X \to \Y$ be a discrete Conduch\'{e} fibration in $\CatE$. Then the following square is a pullback.
  \begin{equation*}
        \begin{tikzcd}
            X_3 \arrow[d, "m^2"'] \arrow[r, "F_3"] & Y_3 \arrow[d, "m^2"] \\
            X_1 \arrow[r, "F_1"'] & Y_1
        \end{tikzcd}
    \end{equation*}

\end{lem}

\begin{proof}
    It is enough to prove this representably in $\Cat,$ in which case it is an easy exercise in lifting a composable triple in $\Y$ that is in the image  of $F$ to a composable triple in $\X$ by using the discrete Conduch\'{e} property twice. 
\end{proof}

   \begin{prop}\label{lem: pullback stable coequalisers that agree on objects}
      Let $\E$ be a category with pullbacks and pullback-stable coequalisers. Then coequalisers of parallel pairs of internal functors that agree on objects are stable under pullback along discrete Conduch\'{e} fibrations. 
   \end{prop}

   \begin{proof}
       Let $F,G: \A \to \B$ in $\CatE/\Y$ with $F_0 = G_0$ and coequaliser $Q: \B \to \C$ and $f: \X \to \Y$ in $\CatE$. We will show that $f^*(Q) : f^*(\B) \to f^*(\C)$ is the coequaliser of $f^*(F), f^*(G) : f^*(\A) \to f^*(\B)$. First, note that it is clear that $f^*(F)_0 = f^*(G)_0$. Following Proposition~\ref{prop coequaliser when agree on objects}, the coequaliser of $f^*(F)$ and $f^*(G)$ is given by the coequaliser of the following diagram

       \begin{equation*}
    \begin{tikzcd}[column sep = huge]
        f^*_3(L) \arrow[r, shift left = 2, "f^*\left(\widetilde{F}\right)"]\arrow[r, shift right = 2, "f^*\left(\widetilde{G}\right)"'] & f^*_3(B_3) \arrow[r, "f^*(m^2)"] & f^*_1(B_1).
    \end{tikzcd}
\end{equation*}

The discrete Conduch\'{e} condition allows us to rewrite the above diagram as a pullback over $f_1$ rather than a mixture of $f_1$ and $f_3$; by the pullback lemma and Lemma~\ref{lem: conduche} the outside of the following diagrams are pullbacks.

\begin{equation*}
    \begin{tikzcd}
        f_3^*(L)\arrow[dr, phantom, "\lrcorner", very near start]  \arrow[r] \arrow[d] & L \arrow[d] \\
        X_3 \arrow[r, "f_3"] \arrow[d, "m^2"'] \arrow[dr, phantom, "\lrcorner", very near start] & Y_3 \arrow[d, "m^2"] \\
        X_1 \arrow[r, "f_1"] & Y_1 
    \end{tikzcd} \qquad 
    \begin{tikzcd}
        f_3^*(B_3)\arrow[dr, phantom, "\lrcorner", very near start]  \arrow[r] \arrow[d] & B_3 \arrow[d] \\
        X_3 \arrow[r, "f_3"] \arrow[d, "m^2"'] \arrow[dr, phantom, "\lrcorner", very near start] & Y_3 \arrow[d, "m^2"] \\
        X_1 \arrow[r, "f_1"] & Y_1 
    \end{tikzcd}
\end{equation*}

Therefore, the coequaliser of $f^*(F)$ and $f^*(G)$ is given by the coequaliser of the following diagram.
\begin{equation*}
    \begin{tikzcd}[column sep = huge]
        f^*_1(L) \arrow[r, shift left = 2, "f^*\left(\widetilde{F}\right)"]\arrow[r, shift right = 2, "f^*\left(\widetilde{G}\right)"'] & f^*_1(B_3) \arrow[r, "f^*(m^2)"] & f^*_1(B_1).
    \end{tikzcd}
\end{equation*}
Now, since coequalisers are stable under pullback in $\E$, the coequaliser of this diagram is $f_1^*(Q): f_1^*(B_1) \to f_1^*(C_1),$ as required.

   \end{proof}

   Note that the previous proposition is trivial in $\Cat$ since discrete Conduch\'{e} functors are examples of exponentiable functors, and so pullback along them has right adjoint, so preserves all colimits. In fact, for any elementary topos $\E$, \cite[Theorem 2.37]{johnstone1977topos} tells us that Conduch\'{e} fibrations are exponentiable in $\CatE$ and so this result is only notable in the cases in which $\E$ is more general.

   In Lemma~\ref{lem: pullback stable codescent coequalisers implies pullback stable coequalisers in E} we prove the converse of the above result. This in particular tells us that Conduch\'{e} fibrations are not always exponentiable in $\CatE$; when $\E$ does not have pullback-stable coequalisers, Conduch\'{e} fibrations in $\CatE$ cannot be exponentiable.

   In light of Corollary~\ref{Corollary coequifiers}, Proposition~\ref{lem: pullback stable coequalisers that agree on objects} implies coequifiers in $\CatE$ are stable under pullback along discrete Conduch\'{e} fibrations. 

   \begin{cor}\label{cor: pullback stable coequifiers}
       Let $\E$ be an extensive category with pullbacks and pullback-stable coequalisers. Then the $2$-category $\CatE$ has coequifiers which are stable under pullback along discrete Conduch\'{e} fibrations.
   \end{cor}

   \begin{proof}
       Let $f: \X \to \Y$ be a discrete Conduch\'{e} fibration in $\CatE$ and consider a pair of parallel internal natural transformations $\overline{\alpha}, \overline{\beta}: F \Rightarrow G: \A \to \B$ in $\CatE/\Y$. By Corollary~\ref{Corollary coequifiers}, the coequifier of $\overline{\alpha}, \overline{\beta}$ in $\CatE$ can be expressed as a coequaliser of a pair of parallel internal functors that agree on objects $\hat{\alpha}, \hat{\beta}: \mathbf{2}_{\E} \times \A \to \B$. This is in $\CatE/\Y$ via the projection $\mathbf{2}_{\E}\times \A \to \A$. By Proposition~\ref{lem: pullback stable coequalisers that agree on objects}, this coequaliser is stable under pullback along $f$. Moreover, notice that $f^*(\mathbf{2}_{\E} \times \A) \cong \mathbf{2}_{\E} \times f^*(\A)$ and so in pulling back, we obtain the coequaliser of the pair $f^*(\hat{\alpha}), f^*(\hat{\beta}): \mathbf{2}_{\E} \times f^*(\A) \to f^*(\B)$ in $\CatE/\X.$ Again, by Corollary~\ref{Corollary coequifiers}, this corresponds exactly to the coequifier of the pair of parallel internal natural transformations $f^*(\overline{\alpha}), f^*(\overline{\beta}): f^*(F) \Rightarrow f^*(G): f^*(\A) \to f^*(\B)$. 
   \end{proof}

\section{Coequalisers of pairs of arrows out of a discrete category}
\label{Coequalisers of pairs of arrows out of a discrete category}

Throughout this section, we assume that $\E$ is a category with pullbacks and pullback-stable coequalisers in which the forgetful functor $\mathcal{U}: \CatE_1 \to \mathbf{Gph}(\E)$ has a left adjoint denoted $\F: \mathbf{Gph}(\E) \to \CatE_1$. In Section \ref{sec: prelim}, we give examples of categories satisfying these conditions; these include elementary toposes with natural numbers objects and list-arithmetic pretoposes. The goal of this Section is to prove that $\CatE$ has coequalisers of pairs of arrows $F, G: A_{0} \to \B$ where $A_{0}$ is a discrete category. Our proof uses the universal property of the free category on a graph, which we state explicitly in Corollary \ref{corollary free category after f and g have been coequalised on objects}, to follow. Our proof also uses the fact that if we have pullback-stable coequalisers in $\E$, then $\CatE$ has coequifiers of any parallel pair of natural transformations, as recorded in Corollary~\ref{Corollary coequifiers}.

\begin{cor}\label{corollary free category after f and g have been coequalised on objects}
    Let $A_{0}$ be a discrete category internal to $\E$ and let $F, G: A_{0} \to \B$ be a parallel pair of internal functors. Form the coequaliser $k_{0}: B_{0} \to C_{0}$ of the parallel pair $F_{0}, G_{0}: A_{0} \to B_{0}$ in $\E$. Consider the graph $\mathcal{G}:= (B_1, C_0, k_0\cdot d_{0}, k_0 \cdot d_{1})$ internal to $\E$. There is a category $\mathbb{F}(\mathcal{G})$ and a morphism of graphs $\eta_{\mathcal{G}}: \mathcal{G} \to \mathcal{U}\mathbb{F}(\mathcal{G})$ with the property that for any internal category $\mathbb{H}$ and morphism of graphs $h: \mathcal{G} \to \mathcal{U}(\mathbb{H})$ there is a unique internal functor $h': \mathbb{F}(\mathcal{G}) \to \mathbb{H}$ satisfying $\mathcal{U}(h') \cdot \eta_\mathcal{G} = h$.
\end{cor}

\begin{proof}
    The morphism of graphs $\eta_{\mathcal{G}}: \mathcal{G} \to \mathcal{U}\mathbb{F}(\mathcal{G})$ is the component of the unit for the adjunction $\mathbb{F} \dashv \mathcal{U}$ at the graph $\mathcal{G}$. The property stated for $\eta_\mathcal{G}: \mathcal{G} \to \mathcal{U}\mathbb{F}(\mathcal{G})$ is precisely the universal property of the unit.
\end{proof}

\begin{lem}\label{lemma k morphism of graphs}
    There is a morphism of graphs $k: \mathcal{U}(\B) \to \mathcal{U}\mathbb{F}(\mathcal{G})$ defined on vertices by the coequaliser $k_{0}: B_{0} \to C_{0}$ of $F_{0}$ and $G_{0}$, and on edges by the edge-assignment $(\eta_\mathcal{G})_{1}: \mathcal{G}_{1} = B_{1} \to \mathbb{F}(\mathcal{G})_{1}$. 
\end{lem}

\begin{proof}
    Since $\eta_\mathcal{G}: \mathcal{G} \to \mathcal{U}\mathbb{F}(\mathcal{G})$ is a morphism of graphs, we see that for $i \in \{0, 1\}$, the equation displayed below holds. 
    
    \begin{equation}\label{equation k morphism of graphs}
        d_{i}^{\mathbb{F}(\mathcal{G})}\cdot(\eta_{\mathcal{G}})_{1} =(\eta_{\mathcal{G}})_0 \cdot d_i^{\mathcal{G}} = k_{0}\cdot d_{i}^\mathbb{B}
    \end{equation}

    This is because $k_{1}\cdot d_{1}^\mathbb{B}: B_{1} \to C_{0}$ is the source of $\mathcal{G}$ and $k_{0} \cdot d_{1}^\mathbb{B}: B_{1} \to C_{0}$ is the target of $\mathcal{G}$ and $(\eta_G)_0 - 1_{C_0}$ (as proven in Lemma \ref{lem counit equal on objects}). This says precisely that $k: \mathcal{U}(\B) \to \mathcal{U}\mathbb{F}(\mathcal{G})$ is well-defined as a morphism of graphs.
\end{proof}

The morphism of graphs $k: \mathcal{U}(\B) \to \mathcal{U}\mathbb{F}(\mathcal{G})$ of Lemma \ref{lemma k morphism of graphs} will typically not be compatible with identity or composition structure. This is rectified by constructing a coequifier ensuring each of these conditions is satisfied.

\begin{lem}\label{lemma coequify identities}
    There is a parallel pair of natural transformations $\overline{\alpha}, \overline{\beta}: k_{0} \Rightarrow k_{0}: \mathbf{disc}(B_0) \to \mathbb{F}(\mathcal{G})$ as displayed below left, whose component assigning morphisms $\alpha, \beta: B_{0} \to \mathbb{F}(\mathcal{G})_{1}$ are given by ${(\eta_{\mathcal{G}})}_{1} \cdot i^\mathbb{B}$ and $i^{\mathbb{F}(\mathcal{G})} \cdot k_{0}$ respectively, as displayed below right.

\begin{equation*}
\begin{tikzcd}\label{eq:firstcoequifier}
        \disc(B_0) \arrow[rr, shift left = 2, "k_0"name=A, bend left]\arrow[rr, shift right = 2,"k_0"'name=B, bend right] && \mathbb{F}(\mathcal{G})
        \arrow[from=A, to=B, Rightarrow, shift right = 3, "\overline{\alpha}"', shorten = 5]\arrow[from=A, to=B, Rightarrow, shift left = 3, "\overline{\beta}", shorten = 5] &{}
    \end{tikzcd}
\begin{tikzcd}
    &B_1\arrow[rd, "(\eta_\mathcal{G})_{1}"]
    \\
    B_{0}
    \arrow[ru, "i^\mathbb{B}"]
    \arrow[rd, "k_{0}"'] &&\mathbb{F}(\mathcal{G})_{1}
    \\
    &C_{0}\arrow[ru, "i^{\mathbb{F}(\mathcal{G})}"']
\end{tikzcd}
\end{equation*}
    
\end{lem}

\begin{proof}
    As $\disc(B_{0})$ is discrete, it suffices to show that $\overline{\alpha}$ and $\overline{\beta}$ respect sources and targets. For $\alpha$ this follows from sources and targets for identities for the category $\B$, while for $\beta$ this follows from the same axioms for the category $\mathbb{F}(\mathcal{G})$.
\end{proof}

\begin{lem}\label{Lemma coequifying composition}
    Let $p: \mathbb{F}(\mathcal{G}) \to \mathbb{I}$ be the coequifier of $\overline{\alpha}$ and $\overline{\beta}$. There is a parallel pair of natural transformations $\overline{\gamma}, \overline{\delta}: p \cdot k_{2} \cdot m^\B \Rightarrow p \cdot m^{\mathbb{F}(\mathcal{G})}\cdot {\eta_{\mathcal{G}}}_{1}: \disc(B_{2}) \to \mathbb{I}$ as displayed below left, whose component assigning morphisms ${\gamma}, {\delta}: B_{2} \to \mathbb{I}_{1}$ are given by $p_{1} \cdot (\eta_\mathcal{G})_{1} \cdot m^\mathbb{B}$ and $m^{\mathbb{I}} \cdot p_{2} \cdot (\eta_\mathcal{G})_{2}$ respectively, as displayed below right.

\begin{equation*}
\begin{tikzcd}\label{eq:secondcoequifier}
        \disc(B_2) \arrow[rr, shift left = 2, " p \cdot k_{2} \cdot m^\B "name=A, bend left]\arrow[rr, shift right = 2,"p \cdot m^{\mathbb{F}(\mathcal{G})}\cdot {\eta_{\mathcal{G}}}_{1}"'name=B, bend right] && \mathbb{I} \qquad
        \arrow[from=A, to=B, Rightarrow, shift right = 3, "\overline{\delta}"', shorten = 5]\arrow[from=A, to=B, Rightarrow, shift left = 3, "\overline{\gamma}", shorten = 5] &{}
    \end{tikzcd}
    \begin{tikzcd}
    &B_1\arrow[rd, "p_1.(\eta_\mathcal{G})_{1}"]
    \\
    B_{2} \arrow[ru, "m^\mathbb{B}"]\arrow[rd, "p_2.(\eta_\mathcal{G})_{2}"'] &&I_1
    \\
    &I_2 \arrow[ru, "m^\mathbb{I}"']
\end{tikzcd}
\end{equation*}
\end{lem}

\begin{proof}
    The proof is similar to that for Lemma \ref{lemma coequify identities}, now using sources and targets for composition for the category $\B$ to prove that $\gamma$ respects sources and targets, and sources and targets for the category $\mathbb{I}$ to prove that $\delta$ respects sources and targets. 
\end{proof}

\begin{lem}\label{lem:welldefinedinternalfunctor}
    Let $t: \mathbb{I} \to \C$ be the coequifier of the natural transformations $\overline{\gamma}$ and $\overline{\delta}$ of Lemma \ref{Lemma coequifying composition}. The morphism of graphs displayed below is well-defined as an internal functor.

    $$\begin{tikzcd}
      Q: = (  \B \arrow[r, "k"] & \mathbb{F}(\mathcal{G}) \arrow[r, "p"] & \mathbb{I} \arrow[r, "t"] & \mathbb{C})
    \end{tikzcd}$$
\end{lem}

\begin{proof}
   Respect for identities is witnessed by the commutativity of the following diagram, in which the left region commutes by the definition of the coequifier $p: \mathbb{F}(\mathcal{G}) \to \mathbb{I},$ and the other regions commute by functoriality of $p$ and $t$.

\begin{equation*}
    \begin{tikzcd}
        B_0 \arrow[r, "k_0= q_0"] \arrow[dd, "i^{\mathbb{B}}"'] & \mathbb{F}(\mathcal{G})_0 \arrow[d, "i^{\mathbb{F}}"] \arrow[r, "p_0"] & I_0 \arrow[r, "t_0"] \arrow[dd, "i^{\mathbb{I}}"] & C_0 \arrow[dd, "i^{\mathbb{C}}"] & \\
        & \mathbb{F}(\mathcal{G})_1 \arrow[dr, "p_1"] & \\
        B_1 \arrow[r, "k_1 = \eta_{\mathcal{G}_1}"'] & \mathbb{F}(\mathcal{G})_1 \arrow[r, "p_1"'] & I_1 \arrow[r, "t_1"'] & C_1
    \end{tikzcd}
\end{equation*}

Respect for composition is witnessed by the commutativity of the following diagram, in which the region on the left commutes by definition of the coequifier $t: \mathbb{I} \to \mathbb{C}$ and the region on the right commutes by functoriality of $t.$

\begin{equation*}
    \begin{tikzcd}
        B_2 \arrow[r, "k_2 : = \eta_{\mathcal{G}_2}"] \arrow[dd, "m^{\mathbb{B}}"'] & \mathbb{F}(\mathcal{G})_2 \arrow[r, "p_2"] & I_2 \arrow[d, "m^{\mathbb{I}}"] \arrow[r, "t_2"] & C_2 \arrow[dd, "m^{\mathbb{C}}"] \\
        & & I_1 \arrow[dr, "t_1"] & \\
        B_1 \arrow[r, "k_1 := \eta_{\mathcal{G}_1}"'] & \mathbb{F}(\mathcal{G})_1 \arrow[r, "p_1"'] & I_1 \arrow[r, "t_1"'] & C_1
    \end{tikzcd}
\end{equation*} 
\end{proof}

\begin{prop}\label{coequaliser out of discrete category}
     The internal functors $F,G: A_0 \to \B$ in $\CatE$ have a coequaliser given by $Q: \B \to \mathbb{C}$, where this internal functor is defined as in Lemma~\ref{lem:welldefinedinternalfunctor}.
\end{prop}

\begin{proof}
   Given an internal functor $R: \B \to \D$ such that $RF = RG,$ we show that there exists a unique internal functor $S: \mathbb{C} \to \mathbb{D}$ satisfying $SQ =R$.

\begin{equation*}
    \begin{tikzcd}
         A_0 \arrow[r, shift left = 2, "F"] \arrow[r, shift right = 2, "G"'] & \mathbb{B} \arrow[r, "Q"] \arrow[dr, "R"'] & \mathbb{C}\arrow[d, dashed, "S"] \\
         & & \mathbb{D}
    \end{tikzcd}
\end{equation*}

 Define $S_0: C_0 \to D_0$ by the universal property of $k_0$ as the coequaliser on objects. Note that there is a morphism of graphs $W: = (S_0, R_1): \mathcal{G} \to \U \mathbb{D}$ as exhibited by the commutativity of the following diagrams:

\begin{equation*}
    \begin{tikzcd}
        B_1 \arrow[dd, "k_0 \cdot d_0 "'] \arrow[rr, "R_1"] \arrow[dr, "d_0^{\mathbb{B}}"] & & D_1 \arrow[dd, "d_0^{\mathbb{D}}"] \\
        & B_0 \arrow[dl, "Q_0"] \arrow[dr, "R_0"] & \\
        C_0 \arrow[rr, "S_0"'] & & D_0.
    \end{tikzcd}\qquad
       \begin{tikzcd}
        B_1 \arrow[dd, "k_0 \cdot d_1 "'] \arrow[rr, "R_1"] \arrow[dr, "d_1^{\mathbb{B}}"] & & D_1 \arrow[dd, "d_1^{\mathbb{D}}"] \\
        & B_0 \arrow[dl, "Q_0"] \arrow[dr, "R_0"] & \\
        C_0 \arrow[rr, "S_0"'] & & D_0.
    \end{tikzcd}
\end{equation*}

Hence, by the adjunction $\mathbb{F} \dashv \U,$ there exists a unique internal functor $W^{\#}: \mathbb{F}(\mathcal{G}) \to \mathbb{D}$ such that $\U(W^{\#})\eta_{\mathcal{G}} = W.$ The commutativity of the following diagram shows that $W^{\#}$ coequifies the natural transformations in Equation~\ref{eq:firstcoequifier}, which induces a unique functor $Y: \mathbb{I} \to \mathbb{D}.$

\begin{equation*}
    \begin{tikzcd}
        B_0 \arrow[r, "i^{\mathbb{B}}"] \arrow[dr, "R_0"] \arrow[dd, "Q_0"'] & B_1 \arrow[ddr, "R_1"] \arrow[r, "\eta_{\mathcal{G}_1}"] & \mathbb{F}(\mathcal{G})_1 \arrow[dd, "W^{\#}_1"] \\
        & D_0 \arrow[dr, "i^{\mathbb{D}}"'] & \\
        C_0 \arrow[r, "i^{\mathbb{F}}"', hook] \arrow[ur, "S_0"] & \mathbb{F}(\mathcal{G})_1 \arrow[r, "W^{\#}_1"'] & D_1
    \end{tikzcd}
\end{equation*}

The commutativity of the following diagram shows that $Y$ coequifies the natural transformations in Equation~\ref{eq:secondcoequifier}, which induces a unique functor $Z: \mathbb{C} \to \mathbb{D}.$ 

\begin{equation*}
    \begin{tikzcd}
        B_2 \arrow[r, "m^{\mathbb{B}}"] \arrow[dd, "\eta_{\mathcal{G}_2}"'] \arrow[dr, "R_2"] & B_1 \arrow[r, "\eta_{\mathcal{G}_1}"] \arrow[ddrr, "R_1"] & \mathbb{F}(\mathcal{G})_1 \arrow[ddr, "W^{\#}_1"] \arrow[r, "p_1"] & I_1 \arrow[dd, "Y_1"] \\
        & D_2 \arrow[drr, "m^{\mathbb{D}}"'] & &  \\
        \mathbb{F}(\mathcal{G})_2 \arrow[r, "p_2"'] \arrow[ur, "W^{\#}_2"] & I_2 \arrow[u, "Y_2"] \arrow[r, "m^{\mathbb{I}}"'] & I_1 \arrow[r, "Y_1"'] & D_1
    \end{tikzcd}
\end{equation*}

By construction, $ZQ= R \mathbb{B} \to \mathbb{C}$ and $R: \mathbb{C} \to \mathbb{D}$ is the unique such functor that does this, as required.

\end{proof}

\section{Coequalisers of arbitrary pairs}
\label{sec: coeqaulisers finish}

In this section, we put together all the work from previous sections in order to show that for $\E$ an extensive category with pullbacks and pullback-stable coequalisers in which the forgetful functor $\mathcal{U}: \CatE_1 \to \mathbf{Gph}(\E)$ has a left adjoint, the $2$-category $\CatE$ has coequalisers of arbitrary pairs of arrows. Moreover this gives a recipe for how to calculate coequalisers in $\CatE.$ We give a proof of this through Lemma \ref{Lemma coequalisers via discretes}, which is a more general statement about coequalisers in $2$-categories $\K$ for which the inclusion of discrete objects $\disc({\K}) \to \K$ is sufficiently well-behaved. Our previous results allow us to apply this lemma to the $2$-category $\CatE$. 

\begin{lem}\label{Lemma coequalisers via discretes}
    Let $\mathcal{K}$ be a $2$-category for which the inclusion of the full-subcategory of discrete objects $\mathbf{disc}: \mathbf{Disc}(\mathcal{K}) \to \mathcal{K}$ has a left adjoint $(-)_{0}$ with counit $\varepsilon: \mathbf{disc}((-)_{0}) \to 1_\mathcal{K}$ and unit which is given component-wise by identities. Suppose $\mathcal{K}$ has coequalisers of any parallel pair $f, g: A \to B$ for which either of the following conditions hold.
    
    \begin{enumerate}
        \item $f_{0} = g_{0}$, or
        \item $A$ is in the image of $\mathbf{disc}$.
    \end{enumerate}
    
    Then $\mathcal{K}$ has all coequalisers.
\end{lem}

\begin{proof}
    Let \begin{tikzcd}
        A \arrow[rr, shift left = 2, "f"]\arrow[rr, shift right = 2, "g"']&& B
    \end{tikzcd} be a parallel pair. By condition (2), $\mathcal{K}$ has the coequaliser of $f \cdot \varepsilon_{A}$ with $g \cdot \varepsilon_{A}$. Let $q: B \to C$ denote this coequaliser; it has the property that $qf \cdot \epsilon_A = qg \cdot \epsilon_A$. Applying $(-)_0$ to this, and by noting that $A_0 = \disc(A_0)_0$ since the unit has identities as its components and by the triangle identities for the adjunction, it follows that $(\epsilon_A)_0 = 1_{A_0}$, so $(qf)_0 = (qf \cdot \epsilon_A)_0 = (qg \cdot \epsilon_A)_0 = (qg)_0$, so by condition (1), $qf$ and $qg$ have a coequaliser, $p: C \to D.$ We claim that $qp: B \to D$ is the required coequaliser of $f$ and $g.$ Certainly, $qpf=qpg$ as they agree on objects and arrows by construction, so it remains to show the universal property of the coequaliser holds. Given $r: B \to E$ such that $rf=rg$, then $rf \cdot \epsilon_A=rg \cdot \epsilon_A$ and so by the universal property of $C$ as a coequaliser of $f \cdot \epsilon_A$ and $g\cdot \epsilon_A$ we get an induced unique arrow $t: C \to E$. But then $t(qf) = rf = rg = t(qg)$ so by the universal property of  $D$ as the coequaliser of $qf$ and $qg$, we get an induced unique arrow $w: D \to E$ such that $wpqf=wpqg$, as required.

\end{proof}

We are now able to verify our main result.

\begin{thm}
\label{thm: CatE has coequalisers}
    Let $\E$ be an extensive category with pullbacks and pullback-stable coequalisers in which the forgetful functor $\mathcal{U}: \CatE_1 \to \mathbf{Gph}(\E)$ has a left adjoint. Then the $2$-category $\CatE$ has finite $2$-colimits.
\end{thm}

\begin{proof}
    From the discussion in Section~\ref{coproducts and copowers}, it suffices to show that $\CatE$ has coequalisers. To do this, we verify that Lemma \ref{Lemma coequalisers via discretes} applies to $\mathcal{K}:= \CatE$. The functor $\disc: \mathbf{Disc}(\Cat(\E)) \simeq \E \to \CatE$ has a left adjoint given by $(-)_0: \CatE \to \E$, with $\disc(E)_0 = E$ for any $E \in \E$ \cite[Remark 2.13]{hughes2024elementarytheory2categorysmall}. By Proposition \ref{prop coequaliser when agree on objects}, condition (1) of Lemma \ref{Lemma coequalisers via discretes} holds while by Proposition \ref{coequaliser out of discrete category}, condition (2) of Lemma \ref{Lemma coequalisers via discretes} holds.
\end{proof}

\begin{remark}
    In Section \ref{sec: prelim}, we explore categories satisfying the required conditions to apply this theorem. Examples include when $\E$ is a list-arithmetic pretopos or elementary toposes with natural numbers objects. In either case, parameterised list objects in $\E$ are needed to form free categories on graphs, which are used in the construction of general coequalisers in $\CatE$. However, it is of interest to describe the coequalisers that exist in $\CatE$ when milder assumptions are made on $\E$, such as just exactness properties between limits and colimits. Let $\E$ have finite limits and colimits and suppose moreover that it is extensive and has pullback-stable coequalisers. Consider a parallel pair of internal functors $F, G: \A \to \B$ and let $Q_{0}: B_{0} \to C_{0}$ denote the coequaliser of $F_{0}$ and $G_{0}$. We briefly describe, without proof, what we believe should be a sufficient condition that is weaker than the existence of the free category on the graph $\mathbb{G}:=$\begin{tikzcd}
        B_{1} \arrow[r, shift right = 2,"Q_{0}\cdot d_{0}"']\arrow[r, shift left = 2, "Q_{0}\cdot d_{1}"]& C_{0}
    \end{tikzcd} but under which the coequaliser of $F$ and $G$ still exists in $\CatE$. We describe this explicitly when $\E := \mathbf{FinSet}$ and leave the generalisation to the internal setting to the interested reader. Let $C_{n} \in \mathbf{Gph}(\E)$ denote the cycle of length $n$; this can be built by first constructing the path of length $n$ using the terminal object and coproducts, and then using a coequaliser to identify the source and target of the path. Then the coequaliser of $F, G: \A \to \B$ exists in $\CatE$ if for all $n\in \N$ and any map $C_n \to \mathbb{G}$, the following lifting problem has a solution in $\mathbf{Gph}(\E)$.

    \begin{equation}
    \begin{tikzcd}
        & \mathcal{U}(\B)\arrow[d]
        \\
        C_{n}\arrow[r]\arrow[ru, dashed]&\mathbb{G} 
    \end{tikzcd}
    \end{equation}

    \noindent This is to say that any cycles which appear in the graph produced by taking equivalence classes of objects in $\B$ already exist in the underlying graph of $\B$ itself. This means that the coequaliser of $F\cdot\varepsilon_\A$ and $G\cdot\varepsilon_\A$ can be formed in $\CatE$, without using parameterised list objects in $\E$. We leave detailed verification of this construction under these milder assumptions to future work.
\end{remark}

\section{A characterisation of when $\CatE$ has $2$-colimits}\label{sec: characterisation}

In this section, we prove a converse to Theorem~\ref{thm: CatE has coequalisers}. It is clear that the existence of finite $2$-colimits in $\CatE$ implies the existence of finite colimits in $\E$; this is spelled out in Lemma~\ref{lem: 2-colims implies colims}.  By Lemma~\ref{rem: free graphs are colims}, the existence of finite $2$-colimits in $\CatE$ implies that the functor $\mathcal{U}: \CatE_1 \to \mathbf{Gph}(\E)$ has a left adjoint. We isolate the property of $\CatE$ having codescent coequalisers which are stable under pullback along discrete Conduch\'{e} fibrations as being important as it implies pullback stability of coequalisers in $\E$. This is recorded in Lemma~\ref{lem: pullback stable codescent coequalisers implies pullback stable coequalisers in E}. We prove that the assumptions of $\E$ being an extensive category with pullbacks and pullback-stable coequalisers in which the forgetful functor $\mathcal{U}: \CatE_1 \to \mathbf{Gph}(\E)$ has a left adjoint is equivalent to the $2$-categorical assumptions that $\CatE$ is an extensive $2$-category with finite $2$-colimits, pullbacks and codescent coequalisers which are stable under pullback along discrete Conduch\'{e} fibrations in Theorem~\ref{thm: characterisation}. We state this in purely $2$-categorical terms without reference to the fact that the $2$-category is of the form $\CatE$ in Theorem~\ref{thm: with Bourke} following \cite{bourke2010codescent}. 

\color{black}

\begin{lem}\label{lem: 2-colims implies colims}
    Let $\E$ be a category with pullbacks and suppose $\CatE$ has finite $2$-colimits. Then $\E$ has  finite colimits.
\end{lem}

\begin{proof}
    This is clear from the adjunctions $\mathbf{disc}\dashv (-)_{0}\dashv \mathbf{indisc}$  that colimits in $\E$ can be calculated in $\CatE$ by applying the $\disc$ functor, calculating the colimit and then applying $(-)_0$, since $(-)_0 \circ \disc = \text{id}_{\E}.$
\end{proof}

    The following observation is noted in the case when $\E = \mathbf{Set}$ in (\cite{bourke2010codescent}, Example 2.6).

\begin{lem}\label{rem: free graphs are colims}
    Let $\E$ be a category with pullbacks, and suppose that $\CatE$ has finite $2$-colimits. Then the forgetful functor $\mathcal{U}: \CatE_1 \to \mathbf{Gph}(\E)$ has a left adjoint.
\end{lem}

\begin{proof}
    
Let $\mathcal{G} = ( G_0, G_1, s, t)$ be an internal graph in $\E$. Since $\CatE$ has finite $2$-colimits, we can construct the coinserter of the following diagram in $\CatE$:

    \begin{equation*}
        \begin{tikzcd}
            \disc(G_1) \arrow[r, bend left, "\disc(s)"] \arrow[r, bend right, "\disc(t)"'] & \disc(G_0) \arrow[r, "Q", dashed] & \mathbb{F}(\mathcal{G}).
        \end{tikzcd}
    \end{equation*}

    This universally coinserts a $2$-cell $Q\disc(s) \Rightarrow Q\disc(t)$, which out of a discrete category means that in $\mathbb{F}(\mathcal{G}),$ there is an actual $1$-cell in $\mathbb{F}(\mathcal{G})$ for any arrow in $G_1$, with source and target as desired. The universal property of the coinserter in this situation is exactly the same as the universal property of the free category.
\end{proof}

  \begin{cor}\label{cor: NNO}
      Suppose $\CatE$ is cartesian closed and has finite $2$-colimits. Then $\E$ is cartesian closed and has a natural numbers object.
  \end{cor}

  \begin{proof}
      If $\CatE$ is cartesian closed, then $\E$ is cartesian closed by \cite[Theorem 4.1]{hughes2024elementarytheory2categorysmall}. By Lemma~\ref{rem: free graphs are colims}, the existence of finite $2$-colimits in $\CatE$ implies that we have a left adjoint to $\mathcal{U}: \CatE_1 \to \mathbf{Gph}(\E)$; by construction, this left adjoint restricts when considering one-object categories and graphs to become a left adjoint to $\mathcal{U}: \mathbf{Mon}(\E) \to \E$, so that $\E$ has free monoids. Note here that we are using that one-object internal graphs are simply objects of $\E$. In a cartesian closed category, having free monoids is equivalent to having a natural numbers object by taking the free monoid on the terminal object; a proof of this is given in \cite[p. 190]{johnstone1977topos} and is due to C.J. Mikkelsen. 
  \end{proof}

Next, we prove a converse to Proposition~\ref{lem: pullback stable coequalisers that agree on objects}, showing a correspondence of exactness conditions. The proof requires the following ``two-point suspension" functor $2[-]: \E \to \CatE$.

\begin{define}
    Given $X \in \E$, the internal category $2[X]$ has object of objects $\mathbf{1}+\mathbf{1}'$ and object of arrows $\mathbf{1}+X+\mathbf{1}'$, in which we write $\mathbf{1}'$ to distinguish the two copies of the terminal object. It has identity assigner given by a coproduct of the coprojection maps $i:= (\iota_{\mathbf{1}}+ \iota_{\mathbf{1}'}): \mathbf{1} + \mathbf{1}' \to \mathbf{1} + \mathbf{1}' + X \cong \mathbf{1} + X+ \mathbf{1}'$, and it has source and target maps given by $d_0: = \iota_{\mathbf{1}} + \iota_{1}\cdot! + \iota_{\mathbf{1}'}: \mathbf{1} + X + \mathbf{1}' \to \mathbf{1} + \mathbf{1}'$ and $d_1:= \iota_{\mathbf{1}} + \iota_{1'}\cdot! + \iota_{\mathbf{1}'}: \mathbf{1} + X + \mathbf{1}' \to \mathbf{1} + \mathbf{1}'$. By extensivity, $(\mathbf{1} + X + \mathbf{1}) \times_{\mathbf{1} + \mathbf{1} }(\mathbf{1} +X +\mathbf{1} ) \cong \mathbf{1}  + X + \mathbf{1}$. and we define $m:= \text{id}_{\mathbf{1} +X+\mathbf{1}' }.$ Given a morphism $f: X \to Y$ in $\E$, the internal functor $2[f]: 2[X] \to 2[Y]$ is defined on objects by $2[f]_0 : = \text{id}_{\mathbf{1}+\mathbf{1} '}$ and on morphisms by $2[f]_1:= \mathbf{1} + f + \mathbf{1}$. 
\end{define} 

It is easy to verify that $\mathbf{2}[-]: \E \to \CatE_{1}$ is well-defined as a functor, and that each internal functor $\mathbf{2}[f]: \mathbf{2}[X] \to \mathbf{2}[Y]$ is a discrete Conduch\'{e} fibration.

\begin{lem}\label{lem: 2[-] preserves coequalisers}
    Let $\E$ be an extensive category with pullbacks. Then the functor $2[-]: \E \to \CatE_{1}$ preserves any coequalisers that exist in $\E$. Moreover, these coequalisers are created by the nerve $N: \CatE_{1} \to [\Delta_{\leq 3}^\text{op}, \E]$.
\end{lem}

\begin{proof}
    Consider a parallel pair $f,g: X \to Y$ in $\E.$ Denote the coequaliser in $\E$ by $q: Y \to C$. Denote the coequaliser of the parallel pair $2[f], 2[g] : 2[X] \to 2[Y]$ in $\CatE$ by $p: 2[Y] \to \C$. The functor $(-)_0: \CatE \to \E$ preserves colimits, so $p_0:  \mathbf{1} + \mathbf{1}' \cong \C_0$.  Certainly, $2[q]$ coequalises $2[f]$ and $2[g]$, so this induces a unique arrow $\C \to 2[C]$. Moreover, by extensivity, the coequaliser of $\mathbf{1}+2[f]+\mathbf{1}'$ and $\mathbf{1}+2[g]+\mathbf{1}'$ is $\mathbf{1}+2[q]+\mathbf{1}': \mathbf{1}+2[Y]+\mathbf{1}' \to \mathbf{1}+2[C]+\mathbf{1}'$, so this induces a unique arrow $r_1:\mathbf{1}+2[C]+\mathbf{1}'\to \C_1.$ It is not hard to check that $(p_0, r_1)$ assembles into an internal functor $\mathbf{1}+2[C]+\mathbf{1}' \to \C$ and that this is inverse to $p: \C \to \mathbf{1}+2[C]+\mathbf{1}'$, finishing the proof.
\end{proof}

\begin{lem}\label{lem: 2[-] preserves pbs}
    Let $\E$ be an extensive category. Then the functor $2[-]: \E \to \CatE$ preserves pullbacks.
\end{lem}

\begin{proof}
    This follows representably from $\E = \s$, where it is straightforward to verify.
\end{proof}

   \begin{prop}\label{lem: pullback stable codescent coequalisers implies pullback stable coequalisers in E}
       Let $\E$ be an extensive category with pullbacks and suppose that the $2$-category $\CatE$ has coequalisers of parallel pairs of internal functors that agree on objects. If these coequalisers are stable under pullback along discrete Conduch\'{e} fibrations,  then $\E$ has pullback-stable coequalisers. 
   \end{prop}

   \begin{proof}
       Let $i: A \to B$ and consider the following coequaliser diagram in $\E/B$:

       \begin{equation*}
       \begin{tikzcd}
         X \arrow[r, shift left = 2, "f"] \arrow[r, shift right = 2, "g"'] \arrow[dr, "x"'] & Y \arrow[r, "q"] \arrow[d, "y"]& C \arrow[dl, "c"] \\
         & B &
    \end{tikzcd}
   \end{equation*}

   By Lemma~\ref{lem: 2[-] preserves coequalisers}, this is sent to a coequaliser diagram in $\CatE/2[B]$ under $2[-];$ note that $2[f]_0 = \text{id}_{\mathbf{1}+ \mathbf{1}'} = 2[g]_0$, and so this is a coequaliser of a parallel pair of internal functors that agrees on objects. Therefore, by assumption, it is stable under pullback along $2[i]: 2[A] \to 2[B]$ which is in particular a discrete Conduch\'{e} fibration, and so by Lemma~\ref{lem: 2[-] preserves pbs} we get a coequaliser diagram in $\CatE/2[A]$ given by the following:

\begin{equation}\label{coequaliser 2[-]}
       \begin{tikzcd}
         2\left[i^*X\right] \arrow[r, shift left = 2, "2{[i^*f]}"] \arrow[r, shift right = 2, "2{[i^*g]}"'] \arrow[dr, "2{[i^*x]}"', bend right] & 2{[i^*Y]} \arrow[r, "2{[i^*q]}"] \arrow[d, "2{[i^*y]}"]& 2{[i^*C]} \arrow[dl, "2{[i^*c]}", bend left] \\
         & 2{[A]} &
    \end{tikzcd}
   \end{equation}

Since $2[f]_0 = 2[g]_0,$  by Lemma~\ref{lem: 2[-] preserves coequalisers}, this is a levelwise coequaliser, so the following diagram is a coequaliser in $\E/ (\mathbf{1} + A + \mathbf{1})$:

 \begin{equation*}
       \begin{tikzcd}[column sep = huge]
         \mathbf{1} + i^*X + \mathbf{1}' \arrow[r, shift left = 2, "\mathbf{1} + i^*f + \mathbf{1}'"] \arrow[r, shift right = 2, "\mathbf{1}+ i^*g + \mathbf{1}'"']  & \mathbf{1} + i^*Y + \mathbf{1}' \arrow[r, "\mathbf{1}+ i^*q + \mathbf{1}'"] & \mathbf{1} + i^*C + \mathbf{1}'
    \end{tikzcd}
   \end{equation*}

   Since coequalisers commute with coproducts, it follows that the following is a coequaliser in $\E/A,$ as required.

   \begin{equation*}
       \begin{tikzcd}
         i^*X \arrow[r, shift left = 2, "i^*f"] \arrow[r, shift right = 2, "i^*g"']  & i^*Y \arrow[r, "i^*q"] & i^*C
    \end{tikzcd}
   \end{equation*}

   \end{proof}

   We therefore have the following characterisation of when $\CatE$ has $2$-colimits.

     \begin{thm}\label{thm: characterisation}
       Let $\E$ be a category with pullbacks. Then $\E$ is extensive, has pullback-stable coequalisers, and the forgetful functor $\mathcal{U}: \CatE_1 \to \mathbf{Gph}(\E)$ has a left adjoint if and only if the $2$-category $\CatE$ is extensive, has $2$-colimits, pullbacks and coequalisers of parallel pairs of functors that agree on objects are stable under pullback along discrete Conduch\'{e} fibrations.  
   \end{thm}

   \begin{proof}
       Extensivity of $\E$ being equivalent to extensivity of $\CatE$ is given by \cite[Lemma 5.2]{hughes2024elementarytheory2categorysmall} and $\E$ has pullbacks if and only if $\CatE$ has pullbacks. If $\E$ is a category with pullbacks and pullback-stable coequalisers in which $\mathcal{U}: \CatE_1 \to \mathbf{Gph}(\E)$ has a left adjoint, then Theorem~\ref{thm: CatE has coequalisers} and Proposition~\ref{lem: pullback stable coequalisers that agree on objects} shows that $\CatE$ has finite $2$-colimits and that coequalisers of parallel pairs of functors that agree on objects are stable under pullback along discrete Conduch\'{e} fibrations. Conversely, suppose that $\CatE$ is extensive, has $2$-colimits, pullbacks and coequalisers of parallel pairs of functors that agree on objects are stable under pullback along discrete Conduch\'{e} fibrations. Then Proposition~\ref{lem: pullback stable codescent coequalisers implies pullback stable coequalisers in E} shows that $\E$ has pullback-stable coequalisers and Lemma~\ref{rem: free graphs are colims} shows that $\mathcal{U}: \CatE_1 \to \mathbf{Gph}(\E)$ has a left adjoint. 
   \end{proof}

   \begin{cor}
       There exist Conduch\'{e} fibrations between double categories that are not exponentiable. 
   \end{cor}

   \begin{proof}
       Let $\E=\Cat$. By Theorem~\ref{thm: characterisation}, if all Conduch\'{e} fibrations of double categories (i.e. in $\Cat(\Cat)$) were exponentiable then all coequaliser diagrams in $\Cat(\Cat)$ would be  stable under pullback along discrete conduch\'{e} fibrations, which would imply that $\Cat$ has pullback-stable coequalisers. However, $\Cat$ does not have pullback-stable coequalisers--- a simple counterexample is given on \cite{ncatlab:exponentials}. 
   \end{proof}

In \cite{niefield2020exponentiability}, Conduch\'{e} fibrations of double categories are called \emph{pre-exponentiable} because they are shown to satisfy a lax exponentiability condition.

   \begin{remark}

   Let 

   \begin{equation*}
       \begin{tikzcd}
         \A \arrow[r, shift left = 2, "F"] \arrow[r, shift right = 2, "G"'] & \B \arrow[r, "Q"] & \C
    \end{tikzcd}
   \end{equation*}

   be a coequaliser diagram in $\CatE$ in which $F_0 = G_0$. By Proposition~\ref{prop coequaliser when agree on objects}, the coequalising map $Q:\B \to \C$ is an isomorphism on objects. In $\CatE$, such functors have a special importance--- they are the codescent morphisms \cite{bourke2010codescent}. We use this to phrase Theorem~\ref{thm: characterisation} in a purely $2$-categorical way, without reference to the fact that the $2$-category is of the form $\CatE.$
   \end{remark}

   \begin{define}\label{def: codescent coequalisers}
       Let $\K$ be a $2$-category. A coequaliser diagram
       \begin{equation*}
       \begin{tikzcd}
         A \arrow[r, shift left = 2, "F"] \arrow[r, shift right = 2, "G"'] & B \arrow[r, "Q"] & C
    \end{tikzcd}
   \end{equation*}
   is called a \emph{codescent coequaliser} if $Q: B \to C$ is a codescent morphism in $\K$.
   \end{define}

   We collect the results of this section so far and combine them with Bourke's characterisation of $2$-categories of the form $\K \simeq \CatE$ in the following.
   
   \begin{thm}\label{thm: with Bourke}
       Let $\E$ be an extensive category with pullbacks and pullback-stable coequalisers in which the forgetful functor $\mathcal{U}: \CatE_1 \to \mathbf{Gph}(\E)$ has a left adjoint. Then the $2$-category $\K : = \CatE$ satisfies the conditions listed below. Conversely, if $\mathcal{K}$ satisfies the conditions listed below, then there is a $2$-equivalence $\K \simeq \mathbf{Cat}\left(\E\right)$ where $\E:= \mathbf{Disc}\left(\mathcal{K}\right)$, in which $\E$ is extensive, has pullbacks and pullback-stable coequalisers and the forgetful functor $\mathcal{U}: \CatE_1 \to \mathbf{Gph}(\E)$ has a left adjoint.

    \begin{enumerate}
        \item $\mathcal{K}$ has pullbacks and powers by $\mathbf{2}$.
        \item $\mathcal{K}$ has codescent objects of categories internal to $\mathcal{K}$ whose source and target maps form a two-sided discrete fibration.
        \item Codescent morphisms are effective in $\mathcal{K}$.
        \item Discrete objects in $\mathcal{K}$ are projective, in the sense of Definition 4.13 of \cite{bourke2010codescent}.
        \item For every object $A \in \mathcal{K}$, there is a projective object $P \in \mathcal{K}$ and a codescent morphism $c: P \rightarrow A$.
        \item $\K$ is extensive.
        \item $\K$ has finite $2$-colimits.
        \item Codescent coequalisers in $\K$ are stable under pullback along discrete Conduch\'{e} fibrations. 
    \end{enumerate}
   \end{thm}

   \begin{proof}
       The properties (1)-(5) are the conditions listed in Theorem 4.18 of \cite{bourke2010codescent}, from which we can deduce that $\K \simeq \mathbf{Cat}\left(\E\right)$ where $\E:= \mathbf{Disc}\left(\mathcal{K}\right)$. The properties (6)-(8) allow us to apply Theorem~\ref{thm: characterisation}.
   \end{proof}

\section{Examples and future work}\label{sec: prelim}

In this section, we give examples of extensive categories with pullback-stable coequalisers in which the forgetful functor $\mathcal{U}: \CatE \to \mathbf{Gpd}(\E)$ has a left adjoint. 

\begin{define}[\cite{maietti2010joyal}, Definition 2.4]
Let $\E$ be a category with finite limits. We say that $\E$ has \emph{parametrised list objects} if for any $X \in \E$, there exists an object $L(X) \in \E$ together with morphisms $r^{X}_0: \mathbf{1} \to L(X)$ and $r^X_1: L(X) \times X \to L(X)$ such that for any $b: B \to Y$ and $g: Y \times X \to Y$, there exists a unique $u: B \times L(X) \to Y$ making the following diagram commute:

\begin{equation*}
    \begin{tikzcd}[column sep = huge]
        B \arrow[r, "(1_B{,} r^X_0 \cdot !_B)"] \arrow[dr, "b"', bend right] & B \times L(X) \arrow[d, "u", dashed] & B \times (L(X) \times X) \arrow[l, "1_B \times r^X_1"'] \arrow[d, "(u \times 1_X) \cdot \sigma", dashed] \\
        & Y & Y \times X \arrow[l, "g"]
    \end{tikzcd}
\end{equation*}

in which $\sigma: B \times (L(X) \times X) \to (B \times L(X))\times X$ is the associative isomorphism of the cartesian product.
\end{define}

\begin{remark}
    We note that for any category $\E$ with parametrised list objects, the assignment $X \mapsto L(X)$ extends to a functor $L: \E \to \E;$ on morphisms $f: X \to Y,$ we define $L(f): L(X) \to L(Y)$ by the universal property of the parametrised list objects, taking $B = \mathbf{1}, Y = L(Y), b=r^Y_0$ and $g= r_1^Y \circ (1_{L(Y)}\times f): L(Y)\times X \to L(Y)$ in the above definition. Moreover, there is a multiplication action $\mu_X: L(X) \times L(X) \to L(X)$ defined by the universal property by taking $B = L(X), Y = L(X), b = 1_{L(X)}$ and $g = r_1^X.$ 

    The maps $\mu_X, r_0^X$ furnish $L(X)$ with the structure of a monoid in $(\E, \times, \mathbf{1}).$
\end{remark}

\begin{example}
    Useful intuition is provided by the case $\E = \s$. For any set $X$, $L(X)$ is defined to be the set of words with alphabet $X$, otherwise known as the free monoid generated by $X$. The morphism $r_0^X: \mathbf{1} \to L(X)$ is given by the empty list. The morphism $r_1^X: L(X)\times X \to L(X)$ takes a word $(x_1...x_n)$ and an element $y \in X$ and outputs the word $(x_1...x_ny)$. The morphism $\mu_X: L(X) \times L(X) \to L(X)$ concatenates two words $((x_1...x_n),(y_1...y_m)) \mapsto (x_1,...x_ny_1...y_m)$. The morphism $\nu_X: X \to L(X)$ takes an element $x \in X$ and forms the singleton word $(x) \in L(X).$
\end{example}

\begin{remark}\label{rem: cc means no parametrized}
    Any category with parametrised list objects has a parametrised natural numbers object by taking $X = \mathbf{1}$. We also remark that if $\E$ is cartesian closed, then the existence of parametrised lists objects (resp. a parametrised natural numbers object) is equivalent to the existence of list objects (resp. a natural numbers object) \cite{johnstone2002sketches1}.
\end{remark}

\begin{define}
    A \emph{locos} is a lextensive category $\E$ with parameterised list objects. If $\E$ is also regular, we call it a \emph{regular locos}. If it is exact, we call it a \emph{list-arithmetic pretopos}. If it is additionally locally cartesian closed, we call it an \emph{arithmetic $\Pi$-pretopos.}
   
\end{define}

In particular, we have the following useful properties.

\begin{itemize}
    \item Any locos is extensive and has finite products, so it is distributive \cite{CARBONI1993145}. It also has finite coproducts.
    \item Any list-arithmetic pretopos has coequalisers \cite[\S 3.9]{maietti2010joyal}, and therefore has finite colimits.  
\end{itemize}

We show that our main result, Theorem~\ref{thm: CatE has coequalisers} is satisfied by any locos with finite pullback-stable coequalisers; that is we show that for a locos $\E$, the forgetful functor $\mathcal{U}: \CatE \to \mathbf{Gph}(\E)$ has a left adjoint. This is essentially proven in \cite[Proposition 7.3]{maietti2010joyal}, which shows this result for list-arithmetic pretoposes; the properties of exactness or even regularity are not used for the construction given there; quotients are not used in the proof, which we describe below. 

\subsection{The free internal category on an internal graph}

\label{sec: free cat}

Throughout this section, let $\E$ be a locos, with notation as given in Section~\ref{sec: prelim}. In this section, we recall the free internal category on an internal graph given in Definition 7.2 of \cite{maietti2010joyal}. The description we give is equivalent but uses categorical language to describe the structure rather than the internal type theory of a list-arithmetic pretopos. In Proposition 7.3 of \cite{maietti2010joyal}, it is proven that this forms a left adjoint to the forgetful functor $\mathcal{U}: \CatE_{1} \to \mathbf{Gph}(\E).$

Let $\mathcal{G}= (G_0, G_1, s, t)$. Define $\FG_0 : = G_0$ and $\FG_1$ as the equaliser of the following diagram:

\begin{equation}
\label{equaliser}
    \begin{tikzcd}
        & LG_0 \times G_0 \arrow[dr, "r_1^{G_0}"] & \\
        G_0 \times L(G_1) \times G_0 \arrow[ur, "! \times L(t) \times 1_{G_0}"] \arrow[dr, "1_{G_0} \times L(s) \times !"'] & & L(G_0) \\
        & G_0 \times L(G_0) \arrow[ur, "r_1^{G_0} \cdot \rho"']& 
    \end{tikzcd}
\end{equation}

where $\rho$ denotes the symmetry isomorphism of the cartesian product $\rho: G_0\times L(G_0) \cong L(G_0) \times G_0$ and $!: G_{0} \to \mathbf{1}$ is the unique map to the terminal object. The identity assigner $i: \FG_0 \to \FG_1$ is induced by the universal property of the equaliser, given that $1_{G_0} \times r^{G_1}_0\cdot ! \times 1_{G_0} : G_0 \to G_0 \times L(G_1) \times G_0$ equalises Diagram~\ref{equaliser}. We define $d_1, d_0: \FG_1 \to G_0$ by the following composites:

\begin{equation*}
     d_1: = \left( \begin{tikzcd}
     \FG_1 \arrow[r] & G_0 \times LG_1 \times G_0 \arrow[r, "\pi_0"] & G_0
    \end{tikzcd} \right)
\end{equation*}
\begin{equation*}
    d_0: = \left( \begin{tikzcd}
      \FG_1 \arrow[r] & G_0 \times LG_1 \times G_0 \arrow[r, "\pi_2"] & G_0 
    \end{tikzcd}\right).
\end{equation*}

The following map

\begin{equation*}
    \begin{tikzcd}
        \FG_1\times_{G_0} \FG_1 \arrow[d] \\
        (G_0 \times LG_1 \times G_0) \times_{G_0} (G_0 \times LG_1 \times G_0 \arrow[d, "\cong"] )\\
        G_0 \times LG_1 \times LG_1 \times G_0 \arrow[d, "1_{G_0} \times \mu_{G_1} \times 1_{G_0}"] \\
        G_0 \times LG_1 \times G_0.
    \end{tikzcd}
\end{equation*}

equalises Diagram~\ref{equaliser}. This therefore induces a map $m: \FG_1\times_{G_0} \FG_1 \to \FG_1$.

\begin{define}[7.2 of \cite{maietti2010joyal}]
    Given an internal graph $\mathcal{G} = (G_0, G_1, s, t)$, we define an internal category $\mathbb{F}\mathcal{G} : = (\FG_0, \FG_1, d_1, d_0,i, m)$.    
\end{define}

Moreover, this internal category is the \emph{free} internal category on an internal graph, forming an adjunction as recorded below. The unit of this adjunction $\eta_{\mathcal{G}}$ is defined by $\eta_{\mathcal{G}_0} : = 1_{G_0} : G_0 \to \FG_0$ and  $\eta_{\mathcal{G}_1}: G_1 \to \FG_1$ which is induced by the universal property of the equaliser, given that $(d_1, \nu_{G_1}, d_1): G_1 \to G_0 \times LG_1 \times G_0$ equalises Diagram~\ref{equaliser}. The counit of the adjunction does an internal version of taking a string of composable arrows and composing them.

\begin{thm}[\cite{maietti2010joyal}, Proposition 7.3]
    \label{Theorem free categories on graphs}
    Let $\E$ be a locos. The assignment $\mathcal{G} \mapsto \FG$ provides a left adjoint to the forgetful functor $\U: \CatE_1 \to \mathbf{Gph}(\E).$
\end{thm}

\begin{remark}
    If $\E$ has countable coproducts, then it is not too hard to prove that for a graph $\mathcal{G} : = (G_0, G_1, s, t)$, the object $\FG_1 \cong \Sigma_{n \in \N} G_n$, where for $n>1$, $G_n$ is its object of composable $n$-arrows:

$$G_n : = \underbrace{G_1 \times_{G_0} ... \times_{G_0} G_1}_{n \text{ times}}.$$

In this case, the proof of Theorem~\ref{Theorem free categories on graphs} using the internal type theory of $\E$ corresponds to a proof using the universal property of the coproduct; internal induction becomes external universal property. This proof is categorically elegant, and illuminates that the proof in \cite{maietti2010joyal} does not need regularity or exactness conditions. 

We do not ask for $\E$ to have countable coproducts as this is not an elementary condition, despite the fact that arithmetic $\Pi$-pretoposes with finite colimits which do not have countable coproducts are hard to construct and do not interact well with other toposes--- see, for example, ( \cite{johnstone2002sketches}, D5.1.7). 
\end{remark}

\begin{remark}
    As mentioned, the description we give for the free internal category on an internal graph is different, but equivalent, to the one given by Maietti in \cite{maietti2010joyal}. We choose this description as it does not rely on using the internal language of a list-arithmetic pretopos, and it does not use coproducts which are indeed not needed for the construction of free internal categories on graphs. We briefly describe how to see the equivalence between the different descriptions, although a full proof is left to the interested reader. The key to this proof is in noting that the object of non-empty lists of $G_1$, denoted $L^*(G_1)$ and described in \cite{maietti2010joyal} using the internal language of $\E$, is isomorphic to $L(G_1)\times G_1$; the isomorphism between them is given by the maps $r_1^X: L(G_1) \times G_1 \to L^*(G_1)$ and $(\text{Bck}, \text{Las}): L^*(G_1) \to L(G_1) \times G_1$, where $\text{Las}: L^*(G_1) \to G_1$ internally takes the last element of a non empty list and $\text{Bck}: L^*(G_1) \to L(G_1)$ takes all elements except for the last one. These maps are described inductively using the internal language of $\E$ in (\cite{maietti2010joyal}, Appenix A). One direction of the isomorphism is shown using the universal property of the product and the list object. The other direction is shown using internal induction on list elements, using the internal language of $\E$. The proof then proceeds by using the fact that $L(G_1) \cong \mathbf{1} + G_1 \times L(G_1)$. This is shown in \cite{johnstone2002sketches}. The proof is finished by noticing that the equalising diagrams constructed give the same equaliser. 
\end{remark}

\subsection{Examples of locoi with finite pullback-stable coequalisers}

Below, we record some examples of suitable categories.

\subsubsection{Locally cartesian closed locoi with coequalisers}

Let $\E$ be a locally cartesian closed locos with coequalisers. By Theorem~\ref{Theorem free categories on graphs}, the forgetful functor $\U: \CatE_1 \to \mathbf{Gph}(\E)$ has a left adjoint. By local cartesian closedness, for any $f: X \to Y$ in $\E$, the pullback functor $f^*: \E/Y \to \E/X$ has right adjoint, and so preserves all colimits, in particular coequalisers. Therefore, any example of such categories allows us to apply Theorem~\ref{thm: CatE has coequalisers} and conclude that $\CatE$ has all finite $2$-colimits. Some key examples of interest are given by fixing $A$ a partial combinatory algebra and consider the category of assemblies $\mathbf{Asm}_{A}$ over this. This also holds in the category of \emph{modest} assemblies, $\mathbf{Mod}_{A}$. These examples will be examined in future work by the first named author and Sam Speight. Note that these are not elementary topos, nor are pretoposes as they are not exact but merely regular.

\subsubsection{List-arithmetic pretoposes}

A list-arithmetic pretopos is by definition an exact, extensive category with parametrised list objects. Moreover, it has pullback-stable coequalisers by the following argument, communicated to us by Peter LeFanu Lumsdaine.

\begin{lem}
    Let $\E$ be a list-arithmetic pretopos. Then $\E$ has pullback-stable coequalisers.
\end{lem}

\begin{proof}
    Examine the proof of the existence of coequalisers in \cite[Proposition 3.10]{maietti2010joyal}. This process uses coproducts, list objects, pullbacks, image factorisation and a quotient by an equivalence relation. All of these are preserved by any functor between list-arithmetic pretoposes. Moreover, all of these conditions are local, and so are preserved by slicing. Therefore, these steps are preserved by pullback. 
\end{proof}

A class of examples of list-arithmetic pretoposes are given by the the syntactic category for any univalent universes of dependent type theory that satisfies axiom K and are closed under the empty type, unit type, sum types, dependent sum types, propositional truncations, quotient sets, and parameterised natural numbers type. Conversely, any list-arithmetic pretopos gives a model of Martin-L\"{o}f type theory which satisfies UIP \cite{streicher1993investigations}. Maietti proposes that list-arithmetic pretoposes are a suitable setting to formulate Joyal's arithemtic universes \cite{maietti2010joyal}, and one can formulate many logical (in)completeness theorems internally to them.

\subsubsection{Arithmetic $\Pi$-pretoposes}

An arithmetic $\Pi$-pretopos is a list-arithmetic pretopos which is locally cartesian closed. As a consequence of Theorem 2.5.17 of \cite{johnstone2002sketches}, any locally cartesian closed positive coherent category with natural numbers object has (parametrised) list objects, so we can replace the need of parametrised list-objects with the existence of a natural numbers object in this case. Hence for $\E$ a model of Palmgren's constructive elementary theory of the category of sets \cite{Palmgren2012Constructivist}, which give categorical models of Bishop's constructive set theory, $\CatE$ has finite $2$-colimits. 

\subsubsection{Elementary toposes with natural numbers object}
Any elementary topos with natural numbers object is an example of an arithmetic $\Pi$-pretopos.  Hence elementary toposes with natural numbers object are a suitable setting for our results. This recovers \cite[Corollary 6.10]{johnstone1978algebraic}, which shows that for an elementary topos with natural numbers object, $\CatE$ has coequalisers. Their proof is different from ours, and we generalise their method in Appendix \ref{appendix: alternative proof}.

Hence, any model of the elementary theory of the category of sets \cite{lawvere2005elementary} is a suitable setting for this work too. This is of interest in relation to \cite{hughes2024elementarytheory2categorysmall}; any model of the elementary theory of the $2$-category of small categories has finite $2$-colimits. 

\subsubsection{Grothendieck toposes}

Any Grothendieck topos is an elementary topos and has a natural numbers object given by the constant sheaf on the natural numbers in $\mathbf{Set}.$ However, these examples are already covered by Proposition \ref{prop: CatE locally presentable} by noting that Grothendieck toposes are locally finitely presentable.

\subsubsection{Further extensions}

For any category satisfying the conditions of Theorem~\ref{thm: CatE has coequalisers}, there are a host of model structures on $\CatE$ given by \cite{everaert2005model}; having finite colimits in $\CatE$ is one of the assumptions necessary for their model structures, and so our work gives more examples of when their theorem can be applied. 

This is related to work by the first named author \cite{hughes2024internalmodel}, in which the present work is used to show that for locally cartesian closed $\E$ satisfying the conditions of Theorem~\ref{thm: CatE has coequalisers}, the $(2,1)$-category $\GpdE$ forms a model of Martin-L\"{o}f type theory. 

The axioms studied in the present paper are not satisfied by $\E=\Cat_1$, the $1$-category of small categories, despite the fact that it is a locos and so by Theorem~\ref{Theorem free categories on graphs}, the forgetful functor $\U: \CatE_1 \to \mathbf{Gph}(\E)$ has a left adjoint. The issue lies in the fact that coequalisers in $\Cat$ are not stable under pullback. An example that shows this is is given on \cite{ncatlab:exponentials}. However, $\Cat$ is locally finitely presentable, so we could apply Proposition \ref{prop: CatE locally presentable} and conclude that $\Cat(\Cat)$, the $2$-category of double categories has finite colimits. 

In future work, we look to extend the method of this paper to prove this without using local finite presentability, and more generally for $\E = \Cat(\mathcal{C})_1$ the $1$-category of categories internal to $\mathcal{C}$ a an extensive $1$-category $\E$ with pullbacks and pullback-stable coequalisers in which the forgetful functor $\mathcal{U}: \CatE_1 \to \mathbf{Gph}(\E)$ has a left adjoint to conclude that $\Cat(\CatE)$, the category of double internal categories has finite $2$-colimits. Part of this work proves that $\Cat$ has lax-pullback-stable coequalisers, a result of independent interest. 

In relation to this, we also hope to work on extending this methodology to show: the existence of finite $2$-colimits in $\mathbf{PsCat}(\K)$ for $\K$ a suitable $2$-category after a suggestion by Bryce Clarke; the existence of finite $2$-colimits in $T$-multicategories after a suggestion by Nathanael Arkor; the existence of finite $2$-colimits in the $2$-category of internal models of an essentially algebraic theory after a suggestion by Peter LeFanu Lumsdaine. For the case that $\E=\s$, the latter is shown by a similar argument to Proposition~\ref{prop: CatE locally presentable} as these are equivalent to locally presentable categories, so this would provide a general way of working with internally locally presentable categories. 

\appendix

\section{A proof of associativity in Proposition~\ref{prop coequaliser when agree on objects}}\label{appendix: proof}

We define $C_3$ as the following pullback.

\begin{equation*}
    \begin{tikzcd}
        C_3 \arrow[r, "\pi_{3,0}"] \arrow[d, "\pi_{3,1}"'] \arrow[dr, phantom, "\lrcorner", very near start] & C_2 \arrow[d, "\pi_1"] \\
        C_2 \arrow[r, "\pi_0"] & C_1.
    \end{tikzcd}
\end{equation*}

To show associativity, we must show that the following diagram commutes

\begin{equation}\label{diagram: ass}
    \begin{tikzcd}
        C_3 \arrow[r, "m \times C_1"] \arrow[d, "\sigma"'] & C_2 \arrow[dd, "m"] \\
        C_1 \times_{B_0} C_2 \arrow[d, "C_1 \times m"'] & \\
        C_2 \arrow[r, "m"] & C_1.
    \end{tikzcd}
\end{equation}

Construct $Q_3: B_3 \to C_3$ by the universal property of $C_3$ as a pullback as in the following diagram

\begin{equation*}
    \begin{tikzcd}
        B_3 \arrow[r, "\pi_{3,0}"] \arrow[d, "\pi_{3,1}"'] \arrow[dr, dashed, "Q_3"] & B_2 \arrow[dr, "Q_2"] & \\
        B_2 \arrow[dr, "Q_2"'] & C_3 \arrow[r, "\pi_{3,0}"] \arrow[d, "\pi_{3,1}"'] \arrow[dr, phantom, "\lrcorner", very near start] & C_2 \arrow[d, "\pi_1"] \\
        & C_2 \arrow[r, "\pi_0"'] & C_1
    \end{tikzcd}
\end{equation*}

in which $Q_3$ exists by the commutativity of the following diagram:

\begin{equation*}
     \begin{tikzcd}
        B_3 \arrow[r, "\pi_{3,0}"] \arrow[d, "\pi_{3,1}"'] \arrow[dr, "\lrcorner", very near start, phantom] & B_2 \arrow[dr, "Q_2"] \arrow[d, "\pi_2"] & \\
        B_2 \arrow[dr, "Q_2"']  \arrow[r,"\pi_2"'] & B_1 \arrow[dr, "Q_1"]  & C_2 \arrow[d, "\pi_1"] \\
        & C_2 \arrow[r, "\pi_0"'] & C_1
    \end{tikzcd}
\end{equation*}

Since coequalisers are assumed to be stable under pullback in $\E$, it follows that $Q_3$ coequalises the pair of parallel arrow $L \times_{B_0} L \times_{B_0} L \to B_3$ $$m^2\cdot \widetilde{F} \times_{B_0} m^2\cdot \widetilde{F} \times_{B_0} m^2\cdot \widetilde{F}, \qquad m^2\cdot \widetilde{G} \times_{B_0} m^2\cdot \widetilde{G} \times_{B_0} m^2\cdot \widetilde{G}.$$

Hence we can appeal to the universal property of the coequaliser: to show that Diagram~\ref{diagram: ass} commutes, it is enough to show that the diagram commutes when precomposed with $Q_3.$ This is witnessed by the following diagram. 

\begin{equation*}
\begin{tikzcd}[column sep = large]
    B_3 \arrow[ddr, "\sigma"']  \arrow[ddrrr, "m \times _{B_0} B_1"] \arrow[rrr, "Q_3"] \arrow[dd, "Q_3"'] & & & C_3 \arrow[r, "m \times_{B_0} C_1 "] & C_2 \arrow[dddd, "m"] \\
    &  & & \mathbf{A} & \\
    C_3 \arrow[d, "\sigma"']&  B_1 \times_{B_0} B_2 \arrow[dr, "B_1 \times m"] &  & B_2 \arrow[d, "m"]  \arrow[uur, "Q_2"'] &  \qquad \\
   C_1 \times_{B_0} C_2 \arrow[d, "1 \times_{B_0} m"'] & \mathbf{B} & B_2  \arrow[r, "m"'] \arrow[dll, "Q_2"] & B_1  \arrow[dr, "Q_1"] & \\
    C_2 \arrow[rrrr, "m"'] & & & & C_1
\end{tikzcd} 
\end{equation*}

In the above, the regions labelled $\mathbf{A}$ and $\mathbf{B}$ are shown to commute by appealing to the universal property of $C_2$ as a pullback of $\pi_0, \pi_1: C_2 \to C_1$, and showing that the regions commute after postcomposing with these projections.

The commutativity of the region $\mathbf{A}$ is shown by the following pair of commutative diagrams.

\begin{equation*}
    \begin{tikzcd}[column sep = large]
        B_3 \arrow[rr, "Q_3"] \arrow[ddd, "m \times_{B_0} {B_1}"'] \arrow[dr, "\pi_{3,0}"]& & C_3 \arrow[r, "m \times_{B_0} {C_1}"] \arrow[d, "\pi_{3,0}"] &  C_2   \arrow[ddd, "\pi_0"] \\
        &B_2\arrow[r, "Q_2"] \arrow[d, "m"]  & C_2 \arrow[ddr,"m" ] & \\
        & B_1  \arrow[drr, "Q_1"]&  & \\
        B_2 \arrow[rr, "Q_2"'] \arrow[ur, "\pi_0"]& &  C_2  \arrow[r, "\pi_0"']& C_1 
    \end{tikzcd}
    \end{equation*}
    \begin{equation*}
       \begin{tikzcd}[column sep = large]
        B_3 \arrow[rr, "Q_3"] \arrow[ddd, "m \times_{B_0} {B_1}"'] \arrow[dr, "\pi_{3,1}"]& & C_3 \arrow[r, "m \times_{B_0} {C_1}"] \arrow[d, "\pi_{3,1}"] &  C_2   \arrow[ddd, "\pi_1"] \\
        &B_2\arrow[r, "Q_2"] \arrow[d, "\pi_1"]  & C_2 \arrow[ddr,"\pi_1" ] & \\
        & B_1  \arrow[drr, "Q_1"]&  & \\
        B_2 \arrow[rr, "Q_2"'] \arrow[ur, "\pi_1"]& &  C_2  \arrow[r, "\pi_1"']& C_1 
    \end{tikzcd}
\end{equation*}

The commutativity of the region $\mathbf{B}$ is shown by the following pair of commutative diagrams.

\begin{equation*}
\begin{tikzcd}[column sep = large]
    B_3 \arrow[r, "\sigma"] \arrow[dd, "Q_3"'] \arrow[dr, "\pi_{3,0}"] & B_1 \times_{B_0} B_2 \arrow[r, "{B_1} \times_{B_0} m "] & B_2 \arrow[r, "Q_2"] \arrow[d, "\pi_0"]  & C_2 \arrow[dddd, "\pi_0"] \\
     & B_2 \arrow[r, "\pi_0"]  \arrow[d, "Q_2"] & B_1 \arrow[dddr, "Q_1"] & \\
    C_3 \arrow[dd, "\sigma"'] \arrow[r, "\pi_{3,0}"] & C_2 \arrow[ddrr, "\pi_0"] & & \\
    & & & \\
     C_1 \times_{B_0} C_2 \arrow[rr, "{C_1} \times_{B_0} m"'] &  & C_2 \arrow[r, "\pi_0"'] & C_1 
\end{tikzcd}
\end{equation*}
\begin{equation*}
\begin{tikzcd}[column sep = large]
    B_3 \arrow[r, "\sigma"] \arrow[dd, "Q_3"'] \arrow[dr, "\pi_{3,1}"] & B_1 \times_{B_0} B_2 \arrow[r, "{B_1} \times_{B_0} m "] \arrow[d, "\pi_{B_2}"] & B_2  \arrow[r, "Q_2"] \arrow[d, "\pi_0"] & C_2 \arrow[dddd, "\pi_0"] \\
     & B_2 \arrow[r, "m"] \arrow[d, "Q_2"] & B_1 \arrow[dddr, "Q_1"] & \\
    C_3 \arrow[dd, "\sigma"'] \arrow[r, "\pi_{3,1}"] & C_2 \arrow[ddrr, "m"] & & \\
    & & & \\
     C_1 \times_{B_0} C_2 \arrow[rr, "{C_1} \times_{B_0} m"'] \arrow[uur, "\pi_{C_2}"] &  & C_2 \arrow[r, "\pi_0"'] & C_1 
\end{tikzcd}
\end{equation*}

Putting all the above steps together, we have shown that associativity holds. 

\section{A second proof using monadicity}\label{appendix: alternative proof}

In this section, we give a second proof of our result using the methodology described in \cite[Corollary 6.10]{johnstone1978algebraic}.

Assume that $\E$ is an extensive category with pullback-stable coequalisers in which the forgetful functor $\mathcal{U}: \CatE_1 \to \mathbf{Gph}(\E)$ has a left adjoint. We wish to show that $\CatE$ has coequalisers of reflexive pairs. To do this, it is enough to show that $\mathcal{U}: \CatE_1 \to \mathbf{Gph}(\E)$ is monadic as monadic functors create all colimits that exist in their codomain and are preserved by the monad. The category $\mathbf{Gph}(\E)$ is an $\E$-valued presheaf category on the parallel arrow category, and so has all colimits that $\E$ does. By assumption, this includes coequalisers.

We wish to apply the crude monadicity theorem to the functor $\mathcal{U}: \CatE \to \mathbf{Gph}(\E).$ The conditions for the crude monadicity theorem are the following.

    \begin{enumerate}
        \item $\mathcal{U}$ has a left adjoint.
        \item $\mathcal{U}$ reflects isomorphisms.
        \item $\CatE$ has coequalisers of reflexive pairs.
        \item $\mathcal{U}$ preserves coequalisers of reflexive pairs. 
    \end{enumerate}

By assumption, (1) holds. By definition, it is clear that (2) holds. Therefore, it remains to show  (3) and (4). To this end, we note the following which follows from Proposition~\ref{prop coequaliser when agree on objects}.

\begin{lem}\label{lem coequaliser equal on objects}
   Let $\E$ be a category with pullback-stable coequalisers. Any reflexive pair $F,G: \mathbb{A} \to \mathbb{B}$ of internal functors between internal categories such that $A_0 = B_0$ and $F_0 = G_0 = 1_{A_0}$ has a reflexive coequaliser in $\CatE$. 
\end{lem}

\begin{remark}
    In this case, the coequaliser is simpler. It is enough to run the argument of Proposition~\ref{prop coequaliser when agree on objects} with $L$ there replaced by the coequaliser of $F_1, G_1: A_1 \to B_1$. This is because we are not causing any newly composable arrows since the categories are already equal on objects. 
\end{remark}

The following is evident in the proof of Proposition~\ref{prop coequaliser when agree on objects}.
\begin{cor}
    The forgetful functor $\mathcal{U}: \CatE \to \mathbf{Gph}(\E)$ creates coequalisers of reflexive pairs $F, G: \A \to \B$ such that $A_0 = B_0$ and $F_0 = G_0 = 1_{A_0}$.  
\end{cor}

The argument in \cite[Corollary 6.10]{johnstone1978algebraic} concludes by saying that both Linton's theorem and the crude monadicity theorem only require these kinds of reflexive coequalisers in their argument. In order to be vaguely self-contained, we spell out the details of this below. 

Since we are in an abstract setting, we must prove the following.

\begin{lem}\label{lem counit equal on objects}
    Let $\E$ be a category with pullback-stable coequalisers and terminal object in which the forgetful functor $\mathcal{U}: \CatE \to \mathbf{Gph}(\E)$ has a left adjoint $\mathbb{F}: \mathbf{Gph}(\E) \to \CatE$. Then for any $\X \in \CatE$, $(\mathbb{F}\mathcal{U} \X)_0 = X_0$.
\end{lem}

   \begin{proof}
       We first note that the functor $(-)_0: \mathbf{Gph}(\E) \to \CatE$ has right adjoint $\codisc: \E \to \mathbf{Gph}(\E)$ which is equal to the composite $\codisc \circ \mathcal{U}: \E \to \CatE \to \mathbf{Gph}(\E).$ Since left adjoints are unique, it follows that $(-)_0 \circ \mathbb{F} = (-)_0$ as required.   
   \end{proof}

We will denote the counit of the adjunction $\epsilon: \mathbb{F}\mathcal{U} \Rightarrow \text{id}$.

Next we note that given $\A \in \CatE$, the following is a reflexive coequaliser diagram that is identity on objects.

\begin{equation*}
    \begin{tikzcd}
        \mathbb{F}\mathcal{U}\mathbb{F}\mathcal{U} \A \arrow[r, shift left = 2, "\epsilon_{\mathbb{F}\mathcal{U}\A}"] \arrow[r, shift right = 2, "\mathbb{F}\mathcal{U}(\epsilon_{\A})"'] & \mathbb{F}\mathcal{U} \A \arrow[r, "\epsilon_{\A}"] & \A. 
    \end{tikzcd}
\end{equation*}

Given $F, G : \A \to \B$, we spell out how to construct its reflexive coequaliser below. Let $Q: \mathcal{U}\mathbb{F} \mathcal{U} \B \to E$ denote the coequaliser in $\mathbf{Gph}(\E)$ of $\mathcal{U}F, \mathcal{U}G: \mathcal{U}\A \to \mathcal{U} \B$ and $Q': \mathcal{U}\B \to E'$ denote the reflexive coequaliser in $\mathbf{Gph}(\E)$ of $$\mathcal{U} \mathbb{F}\mathcal{U}F, \mathcal{U}\mathbb{F}\mathcal{U}G: \mathcal{U}\mathbb{F}\mathcal{U}\A \to \mathcal{U} \mathbb{F}\mathcal{U} \B.$$ 

Note that there is a pair of maps $\F E' \to \F E$ induced by $$\epsilon_{\mathbb{F}\mathcal{U} \B} \mathbb{F}Q : \mathbb{F}\mathcal{U}\mathbb{F}\mathcal{U} \B \to \mathbb{F}E$$

and $$\mathbb{F}\mathcal{U} (\epsilon_{\B}) \mathbb{F}Q: \mathbb{F}\mathcal{U}\mathbb{F}\mathcal{U}\B \to \mathbb{F}E$$

given that they coequalise $\F \mathcal{U} \F \mathcal{U} F$ and $\F \mathcal{U} \F \mathcal{U} G$ by the triangle identities for adjunctions. Denote these induced maps by $u, v: \F E' \to \F E$. Note that these maps are identity-on-objects since all the other maps involved in their construction are too by Lemma~\ref{lem counit equal on objects}.

Calculate this pair's reflexive coequaliser:

\begin{equation*}
    \begin{tikzcd}
        \mathbb{F}E' \arrow[r, shift left = 2, "u"] \arrow[r, shift right =2, "v"'] &  \mathbb{F}E \arrow[r, dashed] & \C
    \end{tikzcd}
\end{equation*}

which exists by Lemma~\ref{lem coequaliser equal on objects} since it is a parallel pair of functors that are identity-on-objects.

Construct the following diagram.

\begin{equation*}
    \begin{tikzcd}[column sep = huge, row sep = huge]
        \mathbb{F}\mathcal{U}\mathbb{F}\mathcal{U}\A \arrow[d, shift left = 2, "\epsilon_{\mathbb{F}\mathcal{U} \A}"] \arrow[d, shift right = 2, "\mathbb{F}\mathcal{U}(\epsilon_{\A})"'] \arrow[r, "\mathbb{F}\mathcal{U}\mathbb{F}\mathcal{U} F", shift left = 2] \arrow[r, "\mathbb{F}\mathcal{U} \mathbb{F}\mathcal{U} G"', shift right = 2] & \mathbb{F}\mathcal{U}\mathbb{F}\mathcal{U}\B \arrow[r] \arrow[d, shift left = 2, "\epsilon_{\mathbb{F}\mathcal{U} \B}"] \arrow[d, shift right = 2, "\mathbb{F}\mathcal{U}(\epsilon_{\B})"']& \mathbb{F}E'  \arrow[d, shift left = 2] \arrow[d, shift right = 2]\\
        \mathbb{F}\mathcal{U}\A \arrow[d, "\epsilon_{\A}"'] \arrow[r, "\mathbb{F}\mathcal{U}F", shift left = 2] \arrow[r, "\mathbb{F}\mathcal{U} G "', shift right = 2] & \mathbb{F}\mathcal{U}\B \arrow[r] \arrow[d, "\epsilon_{\B}"]& \mathbb{F}E \arrow[d] \\
        \A \arrow[r, shift left =2, "F"] \arrow[r, shift right = 2, "G"'] & \B  \arrow[r, dashed] & \C
    \end{tikzcd}
\end{equation*}

In any such diagram in which the top two rows are reflexive coequalisers, and  every column is a reflexive coequaliser, then there exists an arrow on the bottom row making it a coequaliser diagram \cite[Proposition 2.11]{wolff1974v}. Hence we have calculated the coequaliser and shown the following result.

\begin{prop}
    The functor $\mathcal{U}: \CatE \to \mathbf{Gph}(\E)$ is monadic.
\end{prop}

\begin{cor}
    $\CatE$ has finite $2$-colimits.
\end{cor}

\bibliographystyle{alpha}

\end{document}